\documentclass[3p,times]{elsarticle}

\usepackage{latexsym}

\usepackage{indentfirst}
\usepackage[english]{babel}
\usepackage{cancel}

\usepackage{lipsum}
\usepackage{array}
\usepackage{color} 
\usepackage{tabularx}
\usepackage{graphicx} 
\usepackage{amsmath}
\usepackage{amssymb}
\usepackage{amsfonts}	
\usepackage{moreverb}
\usepackage{dsfont}
\usepackage{tipa}
\usepackage{upgreek}
\usepackage{bm}
\usepackage{multirow}
\usepackage{soul}
\usepackage{ textcomp }
\usepackage[dvipsnames]{xcolor}
\usepackage{mathtools}
\usepackage{setspace}
\usepackage{algorithm}
\usepackage{algorithmic}
\usepackage{caption}
\usepackage{subcaption}
\usepackage{booktabs}
\usepackage{cancel}

\allowdisplaybreaks

\newcommand\BibTeX{{\rmfamily B\kern-.05em \textsc{i\kern-.025em b}\kern-.08em
T\kern-.1667em\lower.7ex\hbox{E}\kern-.125emX}}

\usepackage{hyperref} 
\hypersetup{
    colorlinks=true,                          
    linkcolor=blue, 
    citecolor=red, 
    urlcolor=blue  } 

\newcommand{\ie}{{i.e.,}~}

\newcommand{\R}{\mathbf{R}}
\newcommand{\Q}{\mathbf{Q}}
\renewcommand{\S}{\mathbf{S}}

\newcommand{\F}{\mathbf{F}}

\newcommand{\A}{\mathbf{A}}

\newcommand{\bdm}{\begin{displaymath}}
\newcommand{\edm}{\end{displaymath}}

\newcommand{\bea}{\begin{eqnarray} }
\newcommand{\eea}{\end{eqnarray} }

\newfont{\numerikEleven}{ecrm1000}
\newfont{\numerikTen}{cmss10}
\newfont{\numerikNine}{cmss9}
\newfont{\numerikEight}{cmss8}

\journal{Communications in Computational Physics}

\begin{document} 
\begin{frontmatter}
\title{An Arbitrary-Lagrangian-Eulerian solver for relativistic detonation waves } 
\author[UniTN]{Sara Rinaldi}
\ead{sara-rinaldi-1@unitn.it}

\author[UniTN]{Olindo Zanotti\corref{cor1}}
\ead{olindo.zanotti@unitn.it}

\author[UniTN,SUSTech]{Michael Dumbser}
\ead{michael.dumbser@unitn.it}

\address[UniTN]{Laboratory of Applied Mathematics, DICAM, University of Trento, via Mesiano 77, 38123 Trento, Italy} 

\address[SUSTech]{Department of Mathematics, Southern University of Science and Technology, \\ Xueyuan Avenue 1088, 518055 Shenzhen, Guangdong, China}

\begin{abstract} \color[rgb]{0,0,0}
In this paper we study the dynamics of relativistic detonation waves theoretically and numerically. The reaction is physically accounted for by an extra term in the definition of the total energy density and by an additional equation for the evolution of the mass fraction of the reactant, while leaving formally unmodified the equations of mass and energy-momentum conservation. In this way, the Rankine-Hugoniot relations maintain the same formal structure of the inert version.
For the numerical solution we use a finite volume ALE scheme with a second order TVD space-time reconstruction for the homogeneous part of the system, where the mesh velocity is chosen equal to the shock speed, and we combine it through operator splitting with a locally implicit algorithm for the potentially stiff reaction source terms that arise in the equation of the reactant.
We furthermore propose a particularly efficient algorithm for the conversion from the conserved to the primitive variables, which for the relativistic Euler equations is known to be nontrivial.
Following this approach, we can successfully solve the Zel'dovich-von Neumann-Doering profile of a relativistic detonation wave, up to Lorentz factors of the shock front $\gamma_s\sim 7$. 
Our analysis
allowed us to highlight a new special relativistic effect, which has remained unnoticed so far. While in  Newtonian detonations the Zel'dovich pressure jump decreases monotonically with the mass flux through the shock front, in the relativistic case it shows a minimum and then rises monotonically as a function of the mass flux. This may have interesting physical implications on the amount of energy that can be extracted from a relativistic detonation wave.

\end{abstract}

\end{frontmatter}


%
\section{Introduction} 
\label{sec.introduction}
Reaction fronts are discontinuities in a fluid flow that are described as moving surfaces where a chemical (or nuclear) reaction takes place. 
Pioneering investigations in the context of classical Newtonian mechanics were carried out by Zel'dovich in \cite{Zeldovich1971,Zeldovich1987}, followed by a large number of further analysis for which an extended review has been elaborated in \cite{Berkenbosch1994TheOR,Alford28122024}. 
Depending on the physical conditions, the propagating front can be either a shock wave, 
behind which the fluid is compressed and ignited, or a flame front, which moves subsonically with respect to the ambient medium and 
driven by diffusion of heat. This different 
phenomenology corresponds to \emph{detonations}, \citep{MenikoffDetonations,Wagner1989,Lai2019}, and \emph{deflagrations} \citep{MenikoffDeflagrations}, respectively.

In addition to the classical (Newtonian) behavior, reaction flows are also studied in the relativistic regime, mainly due to their relevance for astrophysics, particularly in the cosmological context and as a possible generation mechanism of a stochastic gravitational wave background
\citep{Steinhardt1982,Rezzolla1996,Caprini2008,Leitao2011,Correia2025}. 
One of the first analyses of detonation waves in special relativistic hydrodynamics can be found in \cite{Cissoko1992},
who considered some general conditions under which the solutions of these equations exist. The relativistic Rankine-Hugoniot relations with combustion were also studied by \cite{Gao2012} and,
more recently, \cite{Harpole2019} addressed the full solution of the Riemann problem for relativistic reaction flows, focusing in particular on potential effects induced by tangential velocities, like those found by \cite{Rezzolla2002}.
When we look at the numerical solution of detonation waves, we find several attempts performed in the Newtonian regime, starting from the first preliminary investigations originated in the military context~\cite{Mader1986,Radun2002}.
Further studies of combustion-driven reaction 
waves, focusing on the interplay between dynamics and chemistry, were 
carried out in~\cite{NIKIFORAKIS1996149}. For instance, \cite{Helzel2000} used a fractional step method with a modified Riemann solver; \cite{Sharpe2000} adopted a  second-order Godunov-type scheme to study pulsating detonations while \cite{KleinDeflagration} propose a hybrid front tracking / front capturing scheme. A nonlinear oscillator model for pulsating detonations was proposed in \cite{KleinDetonation}.

A considerable step forward in accuracy of high order schemes for hyperbolic conservation laws with stiff reaction source terms was obtained by  \cite{HidalgoDumbser}, who used the ADER schemes of \cite{schwartzkopff-dumbser-munz,ADERNSE,AMR3DCL}
to solve nonlinear systems of stiff advection--diffusion--reaction equations.
On the contrary, when we move to the relativistic regime, the only  numerical study of detonations, that  we are aware of, is that of \cite{Harpole2019}. 

The aim of this paper is to provide a detailed analysis of relativistic detonation waves, both semi-analytically, through the solution of the  Zel'dovich-von Neumann-Doering (ZND) profile~\cite{fickett_davis_detonation} of an isolated traveling wave, and numerically, via a truly time dependent solution of the corresponding nonlinear system of conservation laws with reaction source terms.
To this extent, we have found  convenient to adopt 
an Arbitrary-Lagrangian-Eulerian (ALE) scheme, originally proposed by~\cite{Hirt1974,Peery2000,Smith1999}, which allows one to solve the dynamics of combustion in a reference frame comoving with the shock front.  A similar approach,
but in the Newtonian regime, has been followed by \cite{Lopato2023} to study the dynamics of detonation waves. We recall that, due to their flexibility, ALE schemes represent a flourishing field of research, with many relevant results obtained with a variety of different numerical schemes and applications~\cite{DumbserUuriintsetsegZanotti2013,Lagrange2D,Lagrange3D,boscheri2014high,Gaburro2020Arepo,ALEDG}.

The plan of the paper is the following: in Sect.~\ref{sec.GHeqs} we present the governing equations of reaction fronts in special relativity, showing how the ZND profile can be computed as the solution of an ordinary differential equation. Sect.~\ref{sec.schemes} is instead devoted to the presentation of the numerical discretization via an ALE finite volume scheme, which includes the treatment of stiff source terms. Sect.~\ref{sec.tests} contains the (one dimensional) numerical results of our investigation, and Sect.~\ref{sec.conclusions} concludes our work.

We adopt a geometrized set of units, in which the speed of light is set to unity, \ie $c=1$, and we assume $(-, +, +, +)$ as signature of the spacetime metric.
Greek indices  denote spacetime indices ranging from 0 to 3, while Latin indices 
are purely spatial, ranging from 1 to 3. Moreover, we set $m_p/k_B=1$, where $k_B$ is the Boltzmann constant and $m_p$ is the molecular mass.

\section{Relativistic reaction shocks}
\label{sec.GHeqs}
In this section we will first present the fundamentals of relativistic fluid dynamics and combustion shock waves in a one-dimensional space setting.
\subsection{Governing equations}
We assume a flat space-time in Cartesian coordinates, with a metric simply given by
\begin{equation}
\mathrm{d}s^2=\eta_{\mu\nu}dx^\mu dx^\nu=-\mathrm{d}t^2+\mathrm{d}x^2+\mathrm{d}y^2+\mathrm{d}z^2\,.
\end{equation}
Even in a framework where reactions take place, we still limit our attention to perfect fluids, which are described by an energy momentum tensor given by 
\begin{equation}
T^{\mu\nu}=(e+p)\,u^{\,\mu}u^{\nu}+p\,\eta^{\,\mu\nu}\,,
\label{eq:T_matter}
\end{equation}
where $u^\mu$ is
the four-velocity of the fluid, while $e$ and $p$ are the energy density and the pressure of the fluid, respectively. The energy density is composed of three terms, namely
\begin{equation}
	e=\rho + \rho\epsilon + \rho q Z\,,
\end{equation}
where $\epsilon$ is the specific (per unit mass) internal energy, while the term $\rho q  Z$ represents the chemical (or nuclear) energy
released during the reaction with $0<q<1$ being the specific energy of the reaction. 
The dimensionless variable $0 \leq Z \leq 1$ denotes
the mass fraction of the  reactant and it varies from $Z=1$, unreacted (unburnt) gas, to $Z=0$, completely reacted (burnt) gas.
Since the enthalpy density is given by $\rho h = e +p$, it follows that the specific enthalpy is affected by the reaction process as
 \begin{equation}
		\label{phrelation-1}
		h=1+\frac{\Gamma}{\Gamma-1}\frac{p}{\rho} +qZ\,. \\
\end{equation}
The equation of state (EOS) is that of an ideal gas\footnote{In the Newtonian regime, detonations for non-ideal equations of state have been considered by \cite{WANG1985}.}, i.e.
 \begin{equation}
	p=\rho \epsilon(\Gamma-1), \label{peosrel}
\end{equation}
where $\Gamma=4/3$ is the adiabatic index of the gas\footnote{
	In a non-degenerate relativistic fluid as described by \cite{synge_1957_rg}, 
	one can define an adiabatic index as
	$\Gamma=1+p/\rho\epsilon$ and find that $\Gamma\rightarrow 4/3$ 
	in the limit of $k_B T\gg m c^2$.
	Moreover, $\Gamma=4/3$
	satisfies Taub's relativistic inequality $(h-p/\rho)(h-4p/\rho)\geq 1$ for any value of the temperature, while $\Gamma=5/3$ does not~\citep{Mignone2007}.}. 
From now on, we will consider a one dimensional motion for the fluid, such that the  four velocity becomes
$u^\mu=\gamma\left(1,v,0,0\right)$, where $v$ is the velocity with respect to the laboratory frame and $\gamma$ is the corresponding Lorentz factor.
As usual, the equations for the conservation of mass and of energy-momentum, which in the relativistic language
are expressed by $\partial_{\mu} (\rho u^{\mu})=0$ and $\partial_{\mu} (T^{\mu\nu})=0$, can then be rephrased as~\citep{Rezzolla_book:2013}
\begin{subequations}
	\begin{align}
		\label{eq-RHD-1}
			\partial_t(\rho \gamma) +\partial_x(\rho \gamma v)&=0\\
		\label{eq-RHD-2}
				\partial_t(h\rho \gamma^2 v)+\partial_x(h\rho \gamma^2 v^2 + p)&=0\\
		\label{eq-RHD-3}
				\partial_t(h\rho \gamma^2 -p)+\partial_x(h\rho \gamma^2 v)&=0.
	\end{align}
\end{subequations}
The conserved variables of this system are:
\begin{equation}
 		D=\rho\gamma, \qquad
 		S=h \rho \gamma^2 v, \qquad
 		E=h\rho \gamma^2 -p.
\end{equation}
In addition to Eqs.~\eqref{eq-RHD-1}-\eqref{eq-RHD-3}, an extra evolution equation is introduced,
governing the variation of the mass fraction of the reactant  as
\begin{equation}
\label{eq-Z}	
	u^\mu\partial_\mu Z = -K(T) Z\,,
\end{equation}
where $K(T)$ is the reaction rate, while
$T$ denotes the temperature, which is computed from the ideal gas EOS  as $T=p/\rho$.
In this paper we consider two different choices for
$K(T)$: 
\begin{itemize}
	\item Discrete ignition temperature model with
	
	\begin{equation}
		\label{Eq:K-ignition}
		K(T)=
		\left\{
		\begin{array}{ll}
			0 \hspace{1cm}\textrm{if}\,\,\, T<T_i\\[10pt]  
			K_0 \hspace{0.8cm}\textrm{if}\,\,\, T\geq T_i\\[10pt]  		
			\end{array}
		\right.
	\end{equation}
	where $T_i$ is the ignition temperature for the combustion.
	\item  Arrhenius law model with
	\begin{equation}
		\label{Eq:K-Arrhenius}
		K(T)=K_0\exp{(-E_a/T)}\,,
	\end{equation} 
	where $E_a$ is called the activation energy.
	
\end{itemize}
Eq.~\eqref{eq-Z} can be written in a fully conservative form, after combining with the continuity equation.
Hence, the complete system of partial differential equations (PDE) to solve is:
 \begin{eqnarray}
 	\label{pdetot-1}
 		\partial_t D+\partial_x(D v)&=&0,\\
 	\label{pdetot-2}
 		\partial_t S+\partial_x(S v+p)&=&0,\\
 	\label{pdetot-3}
 		\partial_t E+\partial_x S&=&0,\\
 	\label{pdetot-4}
 		\partial_t (Z\rho\gamma) + \partial_x(Z \rho \gamma v)&=& -K(T)\rho Z.	
 \end{eqnarray}
A word of caution is also necessary at this stage: the model presented so far should be regarded as a first step towards  a more realistic   description of  detonation waves, whose dynamics is likely to be affected by radiation processes~\cite{Coulombel2012}, which are out of scope of the study conducted in the present paper.

\subsection{Rankine-Hugoniot conditions}
In this section we review  the Rankine-Hugoniot conditions for the PDE system \eqref{pdetot-1}--\eqref{pdetot-4} (see also \cite{Gao2012}). 
As it is customary in this context, we consider a reference frame comoving with the shock front.
Let $V_{s}>0$ be the constant speed of the shock front propagating to the right, and 
$w$ the speed of the fluid with respect to it. We therefore define
\begin{equation}
	w=\frac{v-V_{s}}{1-V_{s} v}<0, \qquad 	 \gamma_w=\frac{1}{\sqrt{1-w^2}}, \qquad v=\frac{w+V_{s}}{1+V_{s} w}. \label{vwequations}
\end{equation}

\begin{figure}[!htbp]
	\begin{center}
		\includegraphics[width=0.6\textwidth, keepaspectratio]{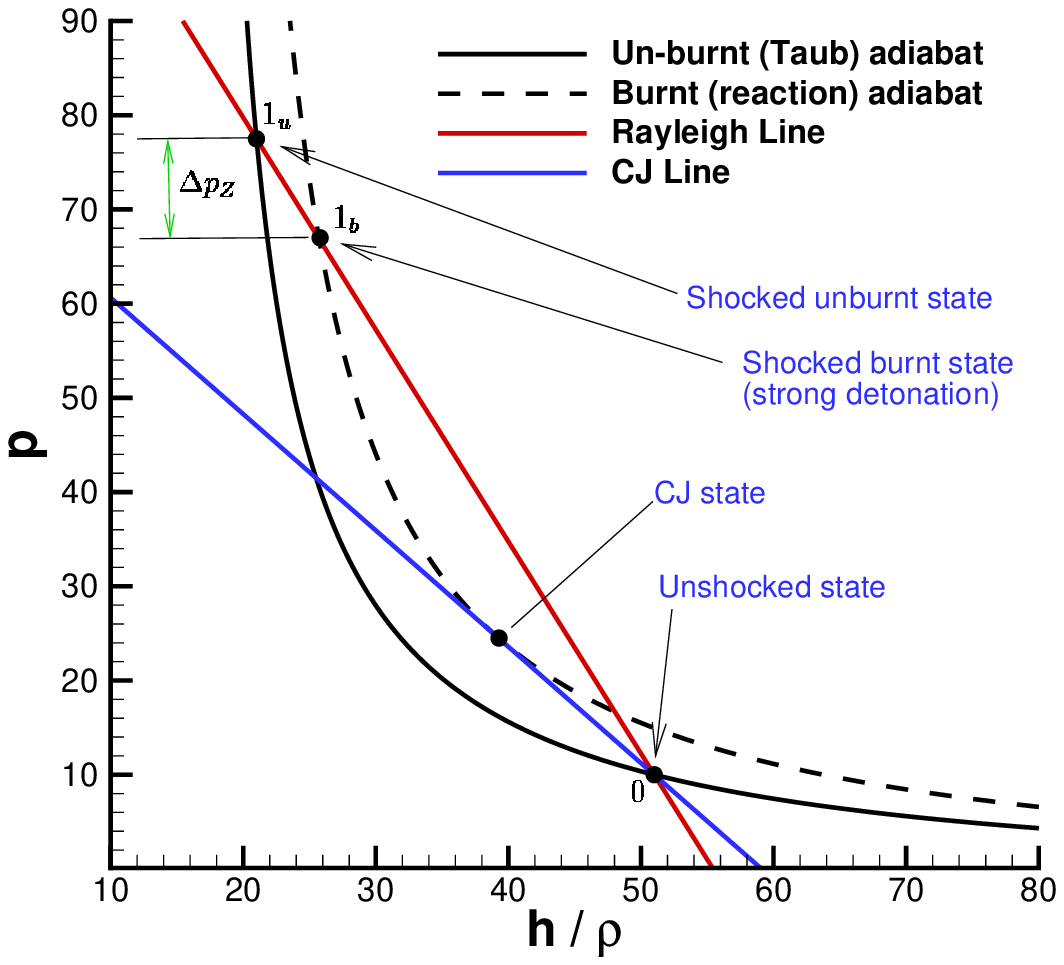}
		\caption{
			Combustion process in the $(h/\rho,p)$ plane for some illustrative values.
		}
		\label{fig:adiabat}
	\end{center}	
\end{figure}

Since the first three equations of the system \eqref{pdetot-1}--\eqref{pdetot-4}
do not have any source term, the standard Rankine-Hugoniot conditions for unreactive flows apply unmodified (see \cite{Rezzolla_book:2013}, Sect. 4.4.3), i.e.
\begin{subequations}
	\begin{align}
		\llbracket \rho \gamma_w w  \rrbracket&=0, \label{RH1}\\ 
		\llbracket h\rho \gamma_w^2  w^2+p \rrbracket&=0, \label{RH2}\\
		\llbracket h\rho \gamma_w^2  w \rrbracket&=0\,,\label{RH3}
	\end{align}
\end{subequations}
where, as usual,  $\llbracket x \rrbracket=x_1-x_0$ indicates the difference of the variable $x$ across the two sides of the shock\footnote{We adopt the convention to indicate the unshocked state with the subscript $0$.}.
From (\ref{RH1}) the mass flux $m=-\gamma_w \rho w>0$ is constant through the shock front. By substituting the mass flux in (\ref{RH2}) and (\ref{RH3}), they are equivalent to:
\begin{subequations}
	\begin{align}
		\llbracket p\rrbracket&= -m^2	\bigg\llbracket \frac{h}{\rho} \bigg\rrbracket, \label{RH2p1}\\
		\llbracket  h \gamma_w \rrbracket&=0,\label{RH3p1}\\
		\llbracket  h^2 \rrbracket&=\left(\frac{h_0}{\rho_0}+\frac{h_1}{\rho_1}\right)\llbracket p\rrbracket\,.\label{RH4p1}
	\end{align}
\end{subequations}
In the presence of combustion, all relations from  \eqref{RH1} to \eqref{RH4p1} still hold, with the reaction energy affecting the thermodynamics through the specific enthalpy $h$. The whole process
can be represented in the $(h/\rho,p)$ plane, as shown in Fig.~\ref{fig:adiabat}, which is the relativistic analog of the $(1/\rho,p$) plane in the classical (Newtonian) case~\citep{Landau-Lifshitz6,Landau1944}. 
All the points in this plane have an enthalpy that can be written as a function of the pressure after solving the quadratic equation
\begin{equation}
	\label{eq:reaction_adiabat}
	\frac{p+p_0(\Gamma-1)}{\Gamma p}h^2+\frac{\Gamma-1}{\Gamma}(1+qZ)\frac{p-p_0}{p}h-h_0^2-\frac{h_0}{\rho_0}(p-p_0)=0\,,	
\end{equation} 
which follows from Eq.~\eqref{RH4p1} by replacing the mass density $\rho$ via the equation of state \eqref{phrelation-1}. See also Eq.~(4.16) of \cite{Pons2000} for a similar expression without the combustion term $qZ$. 
For each value of $Z$, Eq.~\eqref{eq:reaction_adiabat} generates a curve in the  $(h/\rho,p)$ plane, and in Fig.~\ref{fig:adiabat} the two extreme cases are drawn.
Firstly, the solid line indicates those states 
that can be connected through a shock to the initial ``$0$'' state, but
for which combustion has not been ignited yet, hence with $Z=1$.
They collectively form the so-called inert \emph{Taub adiabat}, after \cite{Taub1978}.  In the figure, a representative of such a state is indicated by ``$1_u$'', with the subscript $u$ specifying that it is an \emph{unburnt} state. 
Secondly, the dashed line indicates those states 
that can be connected through a shock to the initial ``$0$'' state, but
for which complete combustion has occurred, hence with $Z=0$.
In the figure, a representative of such a state is indicated by ``$1_b$'', with the subscript $b$ specifying that it is a completely \emph{burnt} state.
The two states ``$1_u$'' and ``$1_b$'' belong to the same Rayleigh line, which comes from Eq.~\eqref{RH2p1} for a constant mass flux $m$.

In the figure we have also reported the so called Chapman-Jouguet (CJ) line, which is tangent to the fully burnt adiabat. We recall that, both in the Newtonian and in the relativistic regime, the CJ line corresponds to the minimum possible mass flux through the detonation wave~\cite{Steinhardt1982}. The computation of the CJ state, including the value of the mass flux $m_{CJ}$, is postponed to Sect.~\ref{sec:JC}. 

\subsection{The relativistic Zel'dovich-von Neumann-Doering profile}
\label{sec:ZNDprofile}
Ahead of the shock, which moves at constant speed $V_s$, no reaction has occurred yet. Just behind the shock front, the temperature rises abruptly, igniting the reaction. This is the state ``$1_u$''. As the combustion goes on, new shocked states all along the Rayleigh line are formed, up until the final fully burnt shocked state
``$1_b$'' is produced. In this way a Zel'dovich-von Neumann-Doering (ZND) profile, characterized by a typical spike (henceforth denoted as the Zel'dovich spike), is produced.

Let us see how this Zel'dovich profile can be computed in a semi-analytic way.
First of all, from (\ref{RH1}) we find the following additional useful kinematic expressions:
\begin{subequations}
	\begin{align}
		w(\rho,m) &=-\frac{m}{\sqrt{m^2+\rho^2}} \label{wrho}\\
		\gamma_w(\rho,m) &=\frac{\sqrt{m^2+\rho^2}}{\rho}\label{odegamma}\\
		v(\rho,m) &=\frac{V_{s}\sqrt{m^2+\rho^2}-m}{\sqrt{m^2+\rho^2}-V_{s}m} \label{vrho}\,.
	\end{align}
	
\end{subequations} 
We now perform a Lorentz transformation into the rest frame of the shock front, i.e.
\begin{equation}
	\label{eq:lorentz}
	\xi = \gamma_s( x-V_{s} t)\,,  \hspace{1cm}\textrm{where}\hspace{1cm}\gamma_s=\frac{1}{\sqrt{1-V_s^2}}
\end{equation}
and we  reduce Eq.~\eqref{eq-Z} into
an ordinary differential equation (ODE)
\begin{equation}
		\partial_{\xi}Z = \frac{K(T) Z}{\gamma\gamma_s(V_s-v)}=\frac{\rho K(T) Z}{m}
		\,,	
			 \label{eq:odez}
\end{equation}
where we have used $\partial_t=-V_s\gamma_s\partial_\xi$ and $\partial_x=\gamma_s\partial_\xi$ and where the mass flux has been written as
\begin{equation}
	m=-\gamma_w \rho w=\rho \gamma\gamma_s(V_s-v)\,,
\end{equation}
which is positive, so that $Z$ decreases from $Z=1$ at the shock front $\xi=0$ towards the burnt state for $\xi<0$.
 Let $V_0=(\rho_0,v_0,p_0,Z_0)$ be the vector of primitive variables of the unshocked (and unburnt) state with $Z_0=1$.
From (\ref{RH1}) we compute the mass flux  as $m=-\gamma_{w,0} \rho_0 w_0$. 
According to (\ref{RH3p1}), 
 $L= h_0 \gamma_{w,0}$ is  constant through the shock, thus leading to
\begin{equation}
	\label{eq:h-rho}
	h(\rho)=\frac{L}{\gamma_w}=\frac{L \rho}{\sqrt{m^2 +\rho^2}}\,.
\end{equation}
Moreover, from the Rankine-Hugoniot condition (\ref{RH2p1}) we can also write:
\begin{equation}
	p(\rho)=p_0 -m^2 \left( \frac{L }{\sqrt{m^2 +\rho^2}}- \frac{h_0}{\rho_0}\right). \label{eq:prho}
\end{equation}
We now replace the pressure $p$ by the equation of state \eqref{phrelation-1}, to find
\begin{eqnarray}
		\frac{\Gamma-1}{\Gamma}\rho(h-1-qZ) &=&p_0 -m^2 \left( \frac{L }{\sqrt{m^2 +\rho^2}}- \frac{h_0}{\rho_0}\right)\\
		\Longrightarrow \frac{\Gamma-1}{\Gamma}\rho\left(\frac{L \rho}{\sqrt{m^2 +\rho^2}}-1-qZ\right) &=&p_0 -m^2 \left( \frac{L }{\sqrt{m^2 +\rho^2}}- \frac{h_0}{\rho_0}\right).
		\label{rhoZ}
\end{eqnarray}
This leads to a fourth degree polynomial equation of the type 
\begin{equation}
	\label{eq:rho4}
	a\rho^4 +b\rho^3 +c \rho^2 +d \rho +e=0\,,
\end{equation} 
with  coefficients that depend only on the state ahead of the shock and on
the value of $Z$ as follows:
\begin{subequations}
	\begin{align}
\label{coeff-a}		
		a &= \left( \frac{\Gamma - 1}{\Gamma} \right)^{2}
		\left( L^{2} - (1 + qZ)^{2} \right)
		\\[8pt]
		b &= -2 \left( \frac{\Gamma - 1}{\Gamma} \right)
		(1 + qZ)\left( p_{0} + m^{2}\frac{h_{0}}{\rho_{0}} \right)
		\\[8pt]
		c &= -\left( p_{0} + m^{2}\frac{h_{0}}{\rho_{0}} \right)^{2}
		+ \left( \frac{\Gamma - 1}{\Gamma} \right) m^{2}
		\left( 2L^{2}
		- \left( \frac{\Gamma - 1}{\Gamma} \right)(1 + qZ)^{2} \right)
		\\[8pt]
		d &= -2 \left( \frac{\Gamma - 1}{\Gamma} \right)
		(1 + qZ)\left( p_{0} + m^{2}\frac{h_{0}}{\rho_{0}} \right) m^{2}
		\\[8pt]
\label{coeff-e}		
		e &= m^{2}\left( m^{2}L^{2}
		- \left( p_{0} + m^{2}\frac{h_{0}}{\rho_{0}} \right)^{2} \right)\,.
	\end{align}
\end{subequations}
The strategy is now the following: we solve numerically the ODE of Eq.~\eqref{eq:odez} using a standard fourth order Runge-Kutta scheme, where the coefficients \eqref{coeff-a}--\eqref{coeff-e} of the quartic for the computation of $\rho$ on the right hand side of \eqref{eq:odez} depend themselves on $Z$.
For each value of $Z$ that is obtained from that ODE, we solve the quartic \eqref{eq:rho4} to obtain $\rho(Z)$, and, 
subsequently, $p(Z)$ and $v(Z)$ from (\ref{eq:prho}) and (\ref{vrho}).
In this way the whole ZND profile can be computed.

\subsection{The Chapman-Jouguet state}
\label{sec:JC}

Unlike the Newtonian case, a closed form expression for the Chapman-Jouguet state is not available in the relativistic regime. As discussed above when commenting on Fig.~\ref{fig:adiabat}, the mass flux $m_{CJ}$ is obtained from the slope of the line that is tangent to the
Taub adiabat of the fully burnt state and passes through the initial unburnt state. 
The numerical computation of the CJ state that we propose here is based on the strategy
outlined in \cite{Anile_book}. We first use  
Eq.~(\ref{RH4p1}) to build the functional
\begin{equation}
	\mathcal{H}:= h^2 -h_0^2 -\left(\frac{h}{\rho}+\frac{h_0}{\rho_0}\right)(p-p_0).
\end{equation}
The curve $\mathcal{H}=0$ describes the Hugoniot curve in the $\left(\frac{h}{\rho},p\right)$ plane.
Its total differential is given by
\begin{align}
		d{\mathcal{H}}&= 2h d{h} -(p-p_0) d{\left(\frac{h}{\rho}\right)}-\left(\frac{h}{\rho}+\frac{h_0}{\rho_0}\right) d{p},\nonumber \\
		&=2hT d{s} -(p-p_0) d{\left(\frac{h}{\rho}\right)}+\left(\frac{h}{\rho}-\frac{h_0}{\rho_0}\right) d{p}\,,
		\label{dH1}
\end{align}
where we have used the first principle of thermodynamics in the form
\begin{subequations}
	\begin{align}
		\label{dh}
		d{h} &= \frac{1}{\rho} d{p} + Td{s},\\ 
		d{\epsilon} &= T d{s} +\frac{p}{\rho^2} d{\rho} \label{deps}\,.
	\end{align}
\end{subequations}
From these we deduce
\begin{subequations}
\begin{align}
		\label{partials_h}
\left. \frac{\partial h }{\partial p} \right|_{s}&=\frac{1}{\rho}, \qquad \qquad \qquad \qquad \left. \frac{\partial h }{\partial s} \right|_{p}=T,\\
\label{partials_eps}
\left. \frac{\partial \epsilon }{\partial \rho} \right|_{s}&=\frac{p}{\rho^2}=c_v \left. \frac{\partial T }{\partial \rho} \right|_{s}, \qquad \left. \frac{\partial \epsilon }{\partial s} \right|_{p}=T+ \frac{p}{\rho^2} \left. \frac{\partial \rho }{\partial s} \right|_{p}\,.
\end{align}
\end{subequations}
After introducing  $\tau=h/\rho$ for ease of notation, 
the Rankine-Hugoniot condition (\ref{RH2p1}) becomes $	\llbracket p\rrbracket / \big\llbracket \tau \big\rrbracket= -m^2	$,
which, in differential form, translates into $dp=-m^2 d\tau$.
Using $(p,\tau)$ as primary thermodynamic variables, we can write
\begin{equation}
	\label{dsdtau}
	\frac{ds}{d\tau}= \left. \frac{\partial s}{\partial \tau} \right|_{p}+ \left. \frac{\partial s}{\partial p} \right|_{\tau} \frac{dp}{d\tau}\,.
\end{equation}
Moreover, from the differential relation $d\tau=\partial\tau/\partial p|_s dp + \partial\tau/\partial s|_p ds$ it follows that
\begin{equation}
	\label{dsdp}
	\left. \frac{\partial s}{\partial p} \right|_{\tau} = -\bigg( \left. \frac{\partial \tau}{\partial s} \bigg)\right|_{p}^{-1}  \left. \frac{\partial \tau }{\partial p} \right|_{s}\,.
\end{equation}
Since the Chapman-Jouguet state is characterized by $\frac{d\mathcal{H}}{d\tau}=0$ and $dp/d\tau=(p-p_0)/(\tau-\tau_0)$, 
we can write, from Eq.~\eqref{dH1}
\begin{align}
\frac{d\mathcal{H}}{d\tau}=&2hT\bigg( \left. \frac{\partial \tau}{\partial s} \bigg)\right|_{p}^{-1} + \bigg( -2hT  \bigg( \left. \frac{\partial \tau}{\partial s}  \bigg)\right|_{p}^{-1} \left. \frac{\partial \tau }{\partial p} \right|_{s} + (\tau-\tau_0) \bigg)\frac{dp}{d\tau} -(p-p_0)=0\nonumber\\
	\Longrightarrow &2hT\bigg( \left. \frac{\partial \tau}{\partial s} \bigg)\right|_{p}^{-1} \mathbin{-}  2hT  \bigg( \left. \frac{\partial \tau}{\partial s}  \bigg)\right|_{p}^{-1} \left. \frac{\partial \tau }{\partial p} \right|_{s}  \frac{dp}{d\tau}=0\,,
	\label{dHdtau}
\end{align}
where we have used both \eqref{dsdtau} and \eqref{dsdp}.
From (\ref{partials_h}) it follows  immediately that:
\begin{align}
		\label{dtaudp}
		\left. \frac{\partial \tau }{\partial p} \right|_{s} &= \frac{1}{\rho}\left. \frac{\partial h }{\partial p} \right|_{s}-\frac{h}{\rho^2} \left. \frac{\partial \rho }{\partial p} \right|_{s} =\frac{1}{\rho^2}-\frac{h}{\rho^2} \left. \frac{\partial \rho }{\partial p} \right|_{s} \\
		\label{dtauds}
		\left. \frac{\partial \tau}{\partial s} \right|_{p}&= \frac{1}{\rho}	\left. \frac{\partial h}{\partial s} \right|_{p}- \frac{h}{\rho^2}\left. \frac{\partial \rho}{\partial s} \right|_{p}=\frac{p}{\rho^2}- \frac{h}{\rho^2}\left. \frac{\partial \rho}{\partial s} \right|_{p}\,.
\end{align}
So far, no assumption has been made on the equation of state adopted. If we now specialize to an ideal gas equation of state, i.e. Eq.~\eqref{peosrel} in (\ref{partials_eps}) we find that:
\begin{equation}
		 \left. \frac{\partial \rho }{\partial p} \right|_{s} =\frac{\rho}{\Gamma p},\qquad
		\left. \frac{\partial \rho}{\partial s} \right|_{p}=-\frac{\rho(\Gamma-1)}{\Gamma },
\end{equation}
Thus Eq.~\eqref{dtaudp}-\eqref{dtauds} become
\begin{subequations}
	\begin{align}
		\label{dtau-1}
		\left. \frac{\partial \tau }{\partial p} \right|_{s} &= \frac{1}{\rho^2}- \frac{\tau}{\Gamma p},\\
		\label{dtau-2}
		\left. \frac{\partial \tau}{\partial s} \right|_{p}&=\frac{p}{\rho^2}+\frac{\tau(\Gamma-1)}{\Gamma },
	\end{align}
\end{subequations}
Replacing \eqref{dtau-1}--\eqref{dtau-2} into
Eq.~(\ref{dHdtau}) we obtain the relation:
\begin{equation}
	\left.\frac{dp}{d\tau}\right|_{CJ}=-m_{CJ}^2=\frac{ \Gamma p\rho^2 }{p \Gamma - h\rho}. \label{eqCJ}
\end{equation}
For any initial state $(\rho_0,v_0,p_0,Z_0=1)$ we solve numerically a system of three equations in three unknowns, namely: Eq.~(\ref{RH2p1}), Eq.~(\ref{eq:rho4}) and Eq.~(\ref{eqCJ}),
to find $\rho$, $p$ and the mass flux $m$ at the Chapman-Jouguet state.
Wherever  the enthalpy $h$ is needed, this is obtained from (\ref{phrelation-1}), while at the CJ state Z is zero and the constant $L$, that is needed for the coefficients of the quartic equation \eqref{eq:rho4}, follows from Eq.~\eqref{eq:h-rho}, i.e.
\begin{equation}
	L=L(m)=\frac{h_0}{\rho_0}\sqrt{m^2+\rho_0^2}\,.
\end{equation}
\begin{table}[!tb]
	\centering
	\caption{Relevant values of the Chapman-Jouguet state as a function of the parameter $q$. Here $p_0=10$, $v_0=0$, $\rho_0=1$. }
	\label{table-CJ}
	\begin{tabular}{llll}
		\toprule
		$q$ & $m_{CJ}$ & $\rho_{CJ}$ & $p_{CJ}$ \\
		\midrule
		0.001 & $0.698453$ & $1.0060743$ & $10.081408$\\
		0.01  & $0.707657$ & $1.0192469$ & $10.260857$\\
		0.1   & $0.737000$ & $1.0611999$ & $10.860056$\\
		0.5   & $0.790687$ & $1.1377872$ & $12.073316$\\
		0.8   & $0.816557$ & $1.1746128$ & $12.718200$\\
		\bottomrule
	\end{tabular}
\end{table}
Following this procedure, it is possible to compute the variables of the CJ state for any set of the initial states ``$0$''.
Table~\ref{table-CJ} reports the CJ variables in a representative case, for different values of the parameter $q$.

\subsection{A new relativistic effect}
\begin{figure}[!htbp]
	\begin{center}
		\includegraphics[width=0.55\textwidth, keepaspectratio]{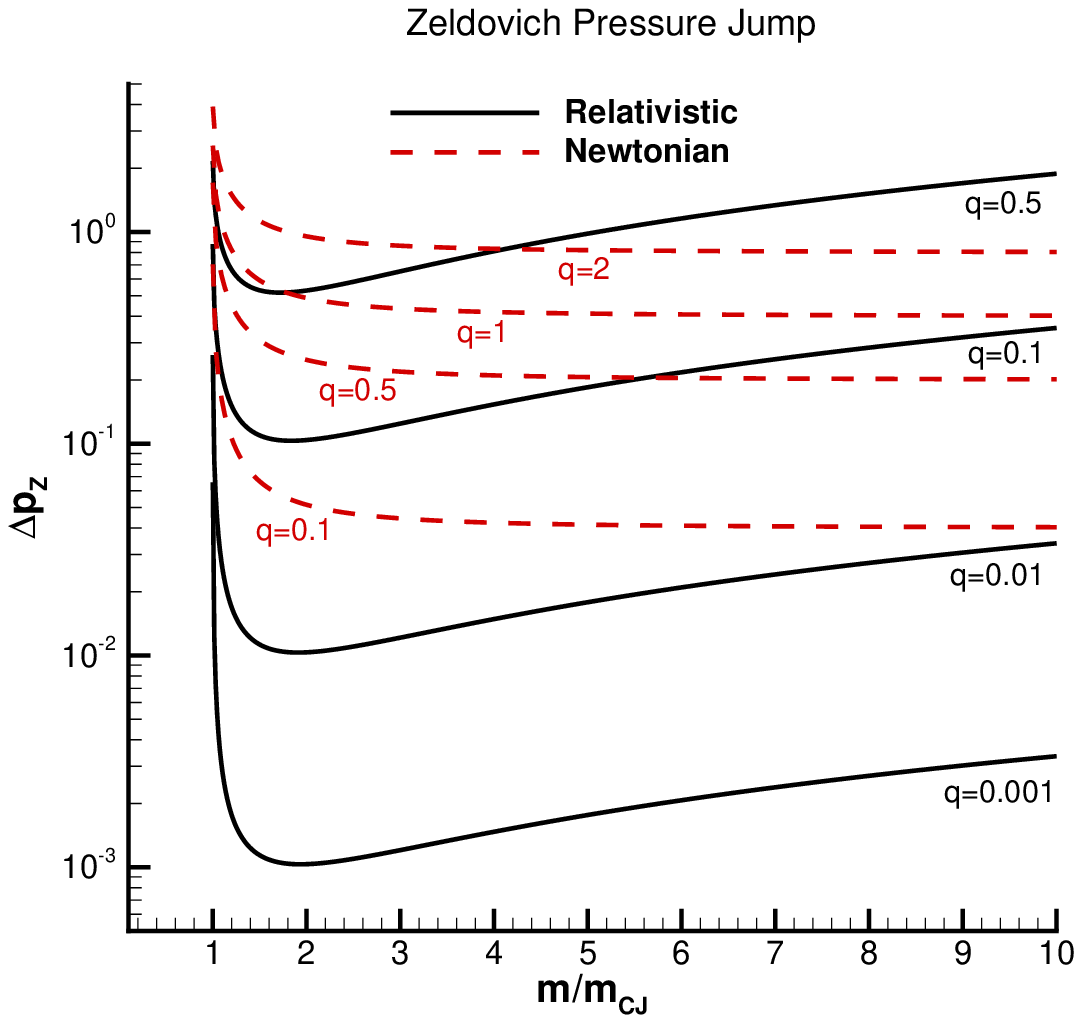}
		\caption{
			Pressure difference between the shocked unburnt and the shocked fully burnt states as a function of the mass flux. Here $p_0=10$, $v_0=0$, $\rho_0=1$.
		}
		\label{fig:Pjump}
	\end{center}	
\end{figure}

\begin{figure}[!htbp]
	\begin{center}
		\includegraphics[width=0.45\textwidth,trim=30 17 20 6,clip]{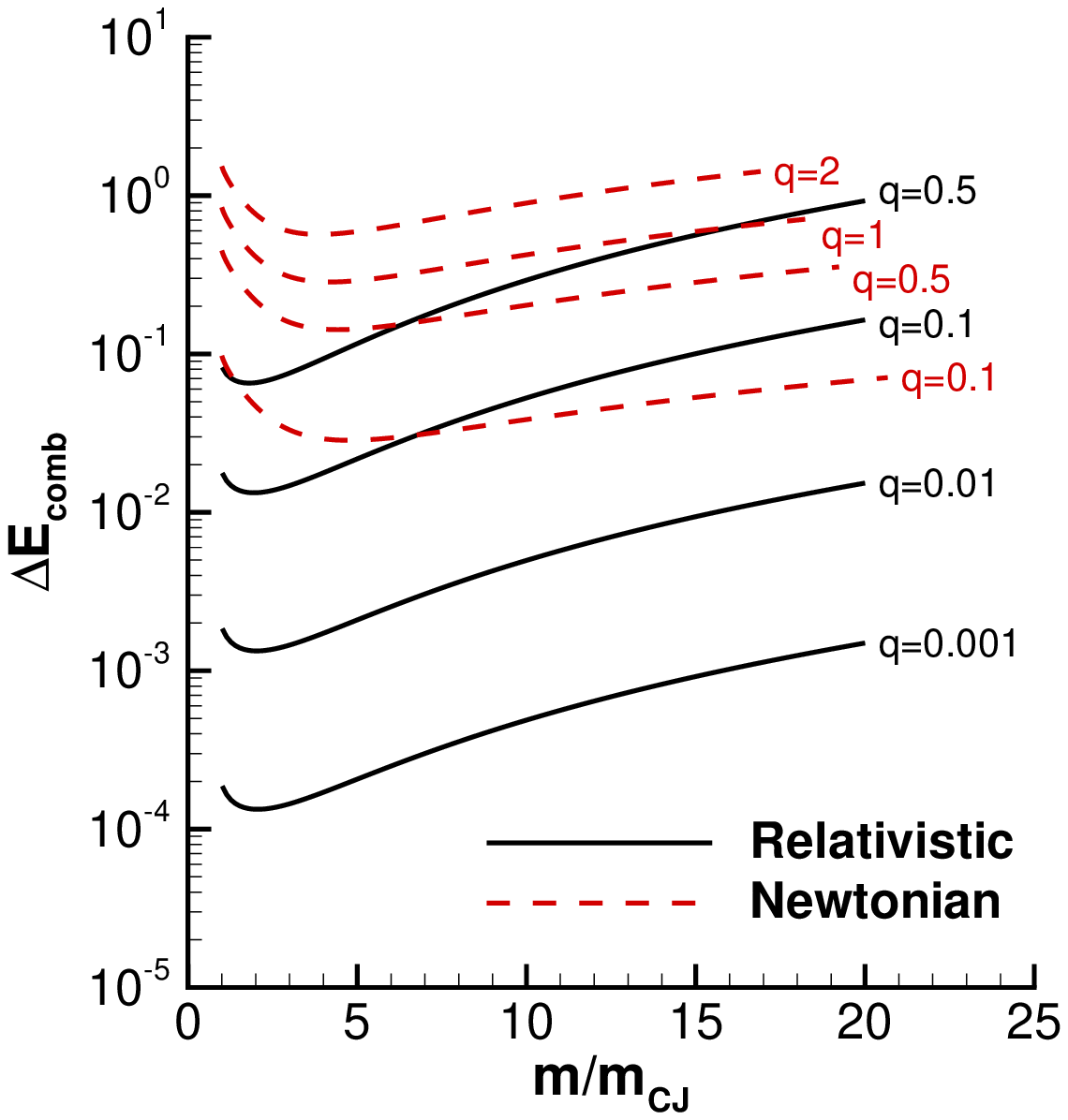}\hfill
		\includegraphics[width=0.45\textwidth,trim=30 17 20 6,clip]{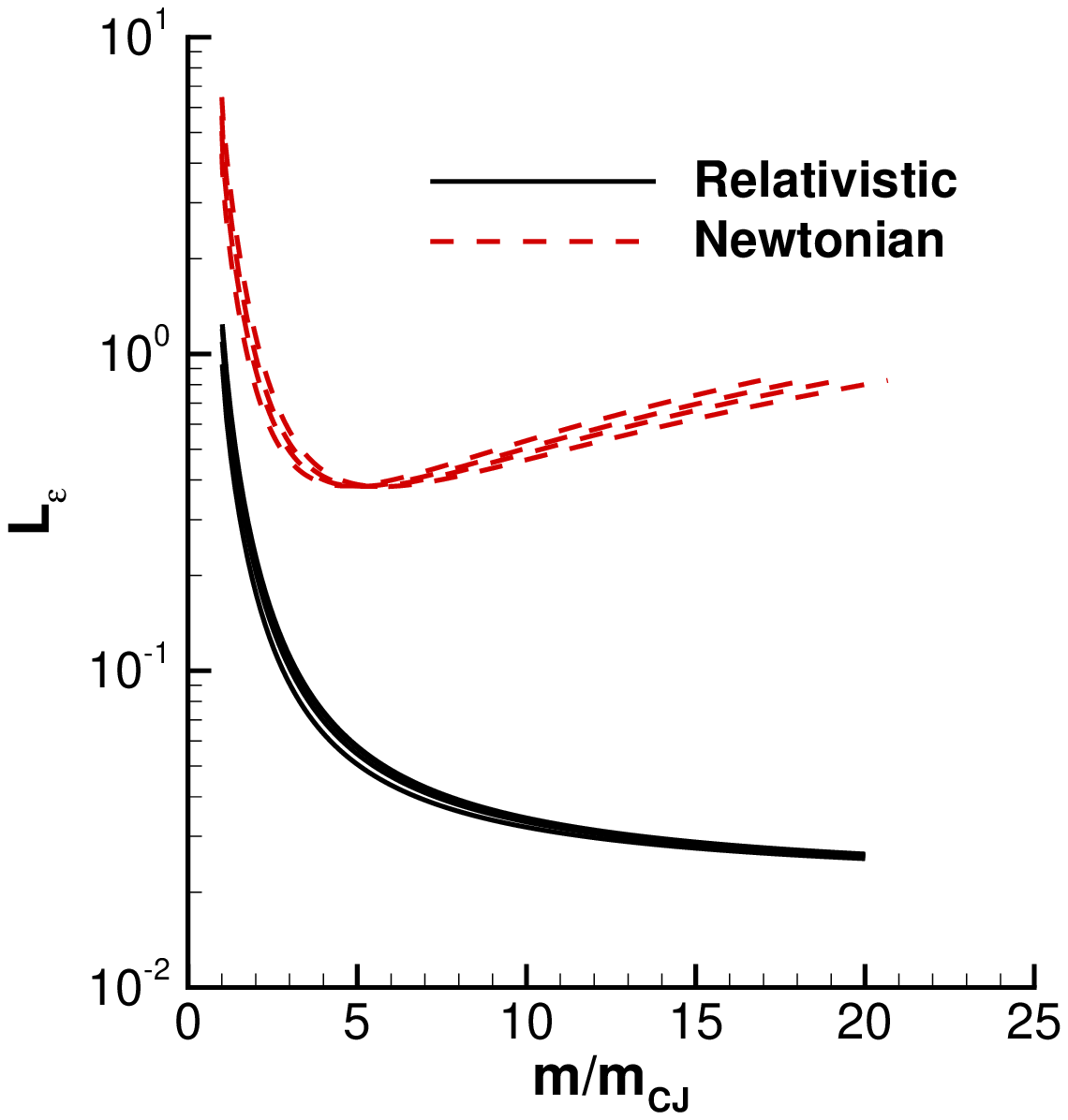}
		\caption{
Chemical energy stored in the ZND profile, computed from \eqref{eq:Ecomb-rel}--\eqref{eq:Ecomb-newt} (left), and width $L_{\varepsilon}$ of the ZND profile for $\varepsilon=10^{-3}$ (right),as a function of the mass . Continuous black lines: relativistic regime; dashed red lines: Newtonian regime. In the right panel the curves for the different values of $q$ nearly coincide. Here $p_0=10$, $v_0=0$, $\rho_0=1$ with an Arrhenius law using $E_a=35$, $K_0=120$.
		}
		\label{fig:EnergyRelease}
	\end{center}	
\end{figure}

Qualitatively, the phenomenology described so far does not differ much from what was
already known in the classical (Newtonian) context. Fig.~4 of \cite{Berkenbosch1994TheOR}, for instance, analyzes classical detonations
reproducing the same physical effects  as those shown in Fig.~\ref{fig:adiabat} above, except for the quantities reported along the axis, which in the relativistic case are  $(h/\rho,p)$ rather than $(1/\rho,p)$. 
However, there is a purely relativistic effect hidden in those figures which is related to the {\emph{Zel'dovich pressure jump}} $\Delta p_Z$. This quantity is  highlighted in Fig.~\ref{fig:adiabat} with a green arrowed vertical line. It represents the pressure difference between the shocked unburnt and the shocked fully burnt states connected by the Rayleigh line. Now, the following crucial difference is observed:
\begin{itemize}
	\item {\emph{Newtonian regime}}: the Zel'dovich pressure jump is maximum at the Chapman-Jouguet state, which, we recall, corresponds to the minimum possible mass flux. When the mass flux through the shock front is increased,  the Zel'dovich pressure jump decreases monotonically to the asymptotic value $\Delta p_Z=(\Gamma-1)q\rho_0$.  
	\item {\emph{Relativistic regime}}: the Zel'dovich pressure  jump $\Delta p_Z$ decreases for small  mass fluxes, down to a minimum value, but then increases again to arbitrary large values for very high mass fluxes. Hence, in the extreme relativistic regime, the Zel'dovich pressure jump can become huge.
\end{itemize}
The comparison among the two regimes is reported in Fig.~\ref{fig:Pjump} for a few generic but representative thermodynamic states.
The continuous black lines and the dashed red lines refer to the relativistic and to the Newtonian regimes, respectively. While in both cases the Zel'dovich pressure jump increases with $q$, the dependence on the mass flux is quite different, as just commented. The comparison has been performed after adopting the same values of the unshocked state, i.e. $p_0=10$, $v_0=0$, $\rho_0=1$, but different adiabatic indices, i.e. $\Gamma=4/3$ and $\Gamma=1.4$, for the relativistic and for the Newtonian regimes, respectively.
Since the units of measure are effectively different in the two regimes, absolute numbers are not relevant in this plot, but one should only pay attention to the different functional dependence on the mass flux.
This purely relativistic effect, which has no analog in the Newtonian case, has also an impact on the chemical energy stored in the reaction zone, which can be computed as
\begin{eqnarray}
	\Delta E_{comb,rel} &=&q\cdot M_{ZND}=q\cdot\int_{\xi_L}^{0} \rho\gamma Z d\xi \,,\label{eq:Ecomb-rel}\\
	\Delta E_{comb,Newt}&=&q\cdot M_{ZND}=q\cdot\int_{\xi_L}^{0} \rho Z d\xi\,,\label{eq:Ecomb-newt}
\end{eqnarray}
with $\gamma$ the Lorentz factor of the fluid in the laboratory frame and with $M_{ZND}$ the mass of the reactant that is contained in the ZND profile. Since $Z$ decays exponentially behind the shock front, the integral over the whole half-space $\xi\leq 0$ is finite. For all the cases of Fig.~\ref{fig:EnergyRelease} we choose $\xi_L=-20$ so that the interval $[\xi_L,0]$ contains essentially the entire ZND profile, with $Z(\xi_L)<10^{-9}$. The width of the ZND profile can be estimated directly from Eq.~\eqref{eq:odez}. Since $\rho$ and $T$ vary only weakly across the reaction zone, $Z(\xi)\approx\exp(\xi/\ell)$ for $\xi\leq 0$, with the characteristic length $\ell=m/(\rho K(T))=\gamma_w|w|/K(T)$ ($\gamma_w=1$ in the Newtonian case). This is the distance travelled by the fluid in the rest frame of the shock front during one reaction time $1/K(T)$. We define the width of the ZND profile as the distance $L_{\varepsilon}$ behind the shock front at which the fraction $1-\varepsilon$ of the reactant has been consumed, i.e. $Z(-L_{\varepsilon})=\varepsilon$, so that $L_{\varepsilon}\approx\ell\ln(1/\varepsilon)$, i.e. $L_{\varepsilon}\approx 6.9\,\ell$ for $\varepsilon=10^{-3}$. By \eqref{eq:odez} we have $\rho Z\,d\xi=(m/K(T))\,dZ$. Evaluating $\gamma$ and $T$ in the shocked burnt state then gives the closed-form estimates $\Delta E_{comb,rel}\approx q\,\gamma\,m/K(T)$ and $\Delta E_{comb,Newt}\approx q\,m/K(T)$. For the parameters of Fig.~\ref{fig:EnergyRelease} these agree with the numerically integrated values within $3\%$ in the relativistic and within $7.5\%$ in the Newtonian case, and within $1\%$ for $m/m_{CJ}\geq 2.5$.
The result of this comparison is shown in Fig.~\ref{fig:EnergyRelease}. While the Zel'dovich pressure jump of Fig.~\ref{fig:Pjump} does not depend on the reaction rate $K(T)$, the chemical energy stored in the ZND profile is indeed affected by $K(T)$, for  which we have assumed the Arrhenius model, adopting $E_a=35$ and $K_0=120$. 
The left panel of Fig.~\ref{fig:EnergyRelease} shows the chemical energy stored in the ZND profile and the right panel the width $L_{\varepsilon}$ of the ZND profile for $\varepsilon=10^{-3}$. In both regimes the stored energy first decreases with the mass flux, reaches a minimum and then increases again, in agreement with the estimate given above. The width $L_{\varepsilon}$ decreases monotonically towards a finite value in the relativistic regime, while in the Newtonian regime it reaches a minimum and then grows with the mass flux.

\section{The numerical method}
\label{sec.schemes}

\subsection{ALE finite volume schemes}
\label{sec:ALE}
In this section we will briefly recall the strategy to 
solve stiff hyperbolic balance laws through 
Arbitrary Lagrangian-Eulerian one-step TVD finite volume schemes, which are second order accurate for the homogeneous part of the system.
The state vector is the typical one for special relativistic hydrodynamics,
with the addition of the combustion term, i.e.
$\Q=(D,S,E,W)$, with $W=\rho\gamma Z = D \, Z$. Moreover, 
fluxes and sources are defined by 
$\F(\Q)=(Dv, Sv +p, S, W v)$ and  $\S(\Q)=(0,0,0,-K(T)\rho Z)$.
In this way,  
the one dimensional PDE \eqref{pdetot-1}-\eqref{pdetot-4} can be cast
as a nonlinear system of balance laws
\begin{equation}
	\frac{\partial \Q}{\partial t} +\frac{\partial \F(\Q)}{\partial x}  = \S(\Q)\label{balancelaw} 
\end{equation} 
defined within
a spatial domain $\Omega(t)$ that moves and changes in time.
Our aim is to follow the fluid as the shock propagates, hence with  a mesh that moves with the speed of the shock.
The computational domain $\Omega$ is discretized by a set of moving mesh points $x_{i+1/2}$ that move with a general mesh velocity $V_{i+1/2}$, thus:
\begin{equation}
	\frac{d}{dt}x_{i+1/2}=V_{i+1/2}.
\end{equation}
The spatial control volumes are defined at the current time $t^n$  as $T^n=[x_{i-1/2}^n, x_{i+1/2}^n]$ where $x_{i+1/2}^n=x_{i+1/2}(t^n)$.
We obtain the integral formulation for the balance law (\ref{balancelaw}) by integrating it in the moving space-time control volume $[x_{i-1/2}(t), x_{i+1/2}(t)]\times[t^n, t^{n+1}]$ and by applying the Gauss Theorem:
\begin{equation}
		\label{eq:generalALEscheme}
	\Delta x^{n+1}_i \Q^{n+1}_i = \Delta x^n_i \Q^n_i 
	- \Delta t \left( \F^V_{i+\frac{1}{2}} - \F^V_{i-\frac{1}{2}} \right) 
	+ \Delta x_i^n \, \Delta t \, \S_i.
\end{equation}
We denoted the time step $\Delta t =t^{n+1}-t^n$, while the mesh spacing at time $t^n$, $\Delta x^n_i=x_{i+1/2}^n-x_{i-1/2}^n$. As in standard finite volume schemes,
the cell average at time \( t^n \), the fluxes and the sources are respectively given by
\begin{eqnarray}
	\Q^n_i &=& \frac{1}{\Delta x^n_i} \int_{x^n_{i-\frac{1}{2}}}^{x^n_{i+\frac{1}{2}}} \Q(x, t^n)  dx\,,\\
		\F^V_{i+\frac{1}{2}} &=& \frac{1}{\Delta t} \int_{t^n}^{t^{n+1}} \left[ \F(x_{i+\frac{1}{2}}(t), t) - V_{i+\frac{1}{2}}(t) \, \Q(x_{i+\frac{1}{2}}(t), t) \right] dt\,,\\
		\S_i &=& \frac{1}{\Delta x^n_i \, \Delta t} \int_{t^n}^{t^{n+1}} \int_{x_{i-\frac{1}{2}}(t)}^{x_{i+\frac{1}{2}}(t)} \S(\Q(x,t))  dx \, dt\,.
\end{eqnarray}
The exact flux $\F^V_{i+1/2}$ is approximated with a numerical flux $\F^V(\Q^-_h, \Q^+_h)$ where  $\Q^-_h=\Q(x_{i+1/2}^-, t)$ and $\Q^+_h=\Q_h(x_{i+1/2}^+, t)$ and:
\begin{equation}
	\F^V(\Q,V)=\F(\Q)-V\Q \quad \textnormal{and} \quad \A^V(\Q,V)=\frac{\partial \F^V}{\partial \Q}.
\end{equation}
Second-order accuracy in space is achieved through a Total Variation Diminishing (TVD) spatial reconstruction using as $\Q_h^+=	\Q^{n+\frac{1}{2},\,+}_{i+\frac{1}{2}}$ and $\Q_h^-=	\Q^{n+\frac{1}{2},\,-}_{i+\frac{1}{2}}$ the following reconstruction:
\begin{equation}
\Q^{n+\frac{1}{2},\,\mp}_{i\pm\frac{1}{2}}
=
\Q^n_i
\pm \frac{1}{2}\,\Delta \Q_i
+
\frac{1}{2}\, \Delta t \,\partial_t \Q_i,
\end{equation}
where:
\begin{align}
	\partial_t \Q_i
	&\approx
	-
	\frac{
		\F^V\left(\Q^{n,-}_{i+\frac{1}{2}}\right)
		-
		\F^V\left(\Q^{n,+}_{i-\frac{1}{2}}\right)
	}{\Delta x_i^n},
	\\
	\Delta \Q_i
	&=
	\operatorname{minmod}
	\bigl(
	\Q^n_{i+1} - \Q^n_i,\;
	\Q^n_i - \Q^n_{i-1}
	\bigr),
	\\
	\Q^{n,\,\mp}_{i\pm\frac{1}{2}}
	&=
	\Q^n_i
	\pm \frac{1}{2}\,\Delta \Q_i.
\end{align}

In our numerical tests we have computed these fluxes following two different approaches
\begin{enumerate}
	\item 
	Rusanov flux:
	\begin{equation}
		\F^{V}_h(\Q_h^-,\Q_h^+)=\frac{1}{2}\left(\F^V (\Q_h^-,V_{i+\frac{1}{2}})+ \F^V(\Q_h^+,V_{i+\frac{1}{2}})\right) -\frac{1}{2}s_{\max}(\Q_h^+-\Q_h^-)
	\end{equation}
	where $s_{\max}=\max\bigl(\max(|\lambda(\A^V(\Q_h^-,V_{i+\frac{1}{2}}))|),\max(|\lambda(\A^V(\Q_h^+,V_{i+\frac{1}{2}}))|)\bigr)$ is the maximum signal speed,
	\item  the Osher-type flux \cite{OsherUniversal}:
	\begin{equation}
		\begin{aligned}
			\F^V_h(\Q^-_h, \Q^+_h) =\ \frac{1}{2} \left[ \F^V(\Q^-_h, V_{i+\frac{1}{2}}) + \F^V(\Q^+_h, V_{i+\frac{1}{2}}) \right] 
			 - \frac{1}{2} \left( \int_0^1 \left| \A^V(\boldsymbol{\psi}(s), V_{i+\frac{1}{2}}) \right| ds \right)(\Q^+_h - \Q^-_h)\,,
		\end{aligned}
	\end{equation}
	where $\boldsymbol{\psi}(s)$ is a straight line segment path connecting the two states  $\Q^-_h$ and $\Q^+_h$, i.e.
\begin{equation}
	\boldsymbol{\psi} = \boldsymbol{\psi}(\Q_h^-, \Q_h^+, s) =
	\Q_h^- + s \left( \Q_h^+ - \Q_h^- \right),
	\qquad 0 \leq s \leq 1\,.
\end{equation}
The absolute value of a matrix is calculated by: 
	\begin{equation}
		|\A|=\R |\Lambda|\R^{-1}
	\end{equation}
	with $\R$ the matrix of the right-eigenvectors and $|\Lambda|$ the diagonal matrix of the absolute values of the eigenvalues of $\A$.
\end{enumerate}
In our specific case, since the mesh velocity is constant in space and in time, the new position of the mesh point $x_{i+1/2}$ at time $t^{n+1}$ is simply given by $x_{i+1/2}^{n+1}=x_{i+1/2}^{n}+ \Delta t V$.
\subsection{Conversion from conserved to primitive variables}
\label{sec:cons2prim}

In order to compute both the fluxes $\F$ and the sources $\S$, given a state vector of conserved variables $\Q$, it is necessary to recover the primitive variables $\rho, v, p$ and $Z$. Except for $Z$, which can be easily computed as $Z=W/D$, the remaining quantities cannot be explicitly reformulated using conservative variables, a fact that represents a well known obstacle in relativistic hydrodynamics, with a variety of solutions proposed over the years \citep{DelZanna2007,Rezzolla_book:2013,Siegel2018,Ripperda2019}. 
Here we follow the same approach as \cite{Cai2024Newton}, which has been proved to be very robust
in the absence of  reactions.
 First of all,  
 we define  the space of admissible solutions:
\begin{equation}
	\mathcal{G}=\left\{ \Q=(D,S,E,W)\  | \ \rho(\Q)>0,\  p(\Q)>0,\  |v(\Q)|<1,\  Z(\Q)\in[0,1] \right\}.
\end{equation}
If $p(\Q)$ is known then:
\begin{eqnarray}
\label{eq:rec-v}
		v(\Q)    &=& \frac{S}{E + p(\Q)}, \\
\label{eq:rec-rho}
		\rho(\Q) &=& D\,\sqrt{1 - v^2(\Q)} .
\end{eqnarray}
Similarly to \cite{WU2015} we now demonstrate that the admissible set $\mathcal{G}$ is equivalent to the following set
\begin{equation}
	\mathcal{G}_1=\left\{ \Q=(D,S,E,W)\  | \ D>0,\ g(\Q):=E-\sqrt{D^2(1+qZ)^2 +S^2}>0, \ \   W/D\in[0,1] \right\},
\end{equation} 
Let us first verify that $\mathcal{G}\subseteq \mathcal{G}_1$, namely that 
if $\Q\in \mathcal{G}$ then $\Q\in \mathcal{G}_1$.\\
First of all, $D=\rho \gamma = \frac{\rho}{\sqrt{1-v^2}}>0$ and $E=\frac{\rho h}{1-v^2}-p\mathrel{\geq} e>0$ if $\Q\in \mathcal{G}$ . As previously mentioned, $W/D=Z$ so $Z\in[0,1] $ if and only if $W/D\in[0,1]$.
We now only need to verify that $E>\sqrt{D^2(1+qZ)^2 +S^2}$ which is true if and only if  $E^{2} > \left(D^{2} (1+qZ)^2 +S^2\right)$. Now, it follows that
\begin{align}
	E^{2} - \left(D^{2}(1+qZ)^2 +S^2\right)
	&=
	\left(\frac{\rho h}{1 - v^{2}} - p\right)^{2}
	- \frac{\rho^{2}}{1 - v^{2}}(1+qZ)^2
	- \left(\frac{\rho h v}{1 - v^{2}}\right)^{2} \nonumber\\[6pt]
	&=
	\left(\frac{\rho h}{1 - v^{2}}\right)^{2}
	+ p^{2}
	- 2p \frac{\rho h}{1 - v^{2}}
	- \frac{\rho^{2}}{1 - v^{2}}(1+qZ)^2
	- \left(\frac{\rho h v}{1 - v^{2}}\right)^{2} \nonumber\\[6pt]
	&=
	\frac{1}{1 - v^{2}}
	\left[
	(\rho h - p)^{2}
	- \rho^{2}(1+qZ)^2
	- p^{2} v^{2}
	\right] \nonumber \\[6pt]
	&\overset{(\ref{phrelation-1})}{=}
	\frac{1}{1 - v^{2}}
	\left[
	\rho^{2}\left(1+\frac{1}{\Gamma-1}\frac{p}{\rho}+q Z\right)^{2}
	- \rho^{2}(1+qZ)^2
	- p^{2} v^{2}
	\right]\nonumber\\[6pt]
	&=
	\frac{1}{1 - v^{2}}
	\left[\frac{1}{(\Gamma-1)^2}p^2+\frac{2}{\Gamma-1}\rho p(1+qZ)
	- p^{2} v^{2}
	\right]  \nonumber\\[6pt]
	&\overset{|v|<1}{>}
	\frac{1}{1 - v^{2}}
	\left[ p^2 \frac{\Gamma(2-\Gamma)}{(\Gamma-1)^2}+
	\frac{2}{\Gamma-1}\rho p(1+qZ)
	\right]  > 0 
\end{align}
for any $\Q\in\mathcal{G}$ and $\Gamma\in (1,2]$. Thus $ g(\Q)>0$ and $\Q \in \mathcal{G}_1$.\\
Secondly, let us  verify that $\mathcal{G}_1\subseteq \mathcal{G}$, namely that
if $\Q\in \mathcal{G}_1$ then $\Q\in \mathcal{G}$.\\
Consider the function of $p$ defined by
\begin{equation}
	\label{phipprim2cons}
	\Phi(p) := \frac{S^{2}}{E+p}
	+ D(1+qZ) \sqrt{1 - \frac{S^{2}}{(E+p)^{2}}} +\frac{p}{\Gamma -1}
	- E ,
	\qquad p \in [0,+\infty),
\end{equation}
with $\Q\in \mathcal{G}_1$.
Obviously, $\Phi(p) \in C^{1}[0,+\infty)$, and
\begin{equation}
	\Phi'(p)
	=
	-\frac{S^2}{(E+p)^2}\left(1-\frac{D(1+qZ)}{\sqrt{(E+p)^2-S^2}}\right)
	+\frac{1}{\Gamma-1}
	\ge
	1 -\frac{S^2}{(E+p)^2}
	> 0,
\end{equation}
for all $p \in [0,+\infty)$, when $E>\sqrt{D^{2}(1+qZ)^2+S^{2}}$ and $\Gamma \in (1,2]$ since $Z=W/D \in[0,1]$. This means that $\Phi(p)$ is a strictly monotonically increasing function of $p$ on $[0,+\infty)$.
Since
\begin{align}
	\Phi(0)
	&=
	\frac{S^{2}}{E}
	+ D (1+qZ)\sqrt{1 - \frac{S^{2}}{E^{2}}}
	- E \nonumber\\
	&=
	\left(  D (1+qZ) -\sqrt{E^2-S^2} \right)\sqrt{\frac{E^2-S^2}{E^2}}
	< 0,
\end{align}
and
\begin{equation}
	\lim_{p \to +\infty} \Phi(p) = +\infty,
	\qquad
	\text{since }
	\lim_{p \to +\infty} \frac{\Phi(p)}{p}
	=
	\frac{1}{\Gamma-1} 
	> 0 ,
\end{equation}
by the intermediate value theorem and the monotonicity of $\Phi(p)$,
there exists a unique positive solution to the equation $\Phi(p)=0$ which, from \eqref{phipprim2cons}, satisfies:
\begin{equation}
		\label{eq:zerophip}
		D(1+qZ) \sqrt{1 - \frac{S^{2}}{(E+p)^{2}}}=E-\frac{p}{\Gamma -1}- \frac{S^{2}}{E+p}
	\end{equation}
		
 and such $p(\Q)$ is also the solution of:
\begin{equation}
	D(1+qZ)\gamma=(E+p) -\frac{\Gamma}{\Gamma-1}p\gamma^2\,.
\end{equation}
For any $\Q\in\mathcal{G}_1$, let $p(\Q)$ denote the solution
	of $\Phi(p)=0$, i.e. of Eq.~\eqref{eq:zerophip}. It then follows immediately that
\begin{equation}
	 |v(\Q)|
	= \frac{|S|}{E+p(\Q)}
	< \frac{|S|}{E}
	< 1,
	\qquad
	\rho(\Q) = D\sqrt{1-v^{2}(\Q)}> 0\,,
\end{equation} which concludes the proof.\\
Another easily proven property of the set $\mathcal{G}_1$ is that it is convex. Following closely the argument of \cite{WU2015},
we verify    that for any $\Q_1,\ \Q_2 \in \mathcal{G}_1 $ also $ \Q_{\lambda}:=\lambda \Q_1 + (1-\lambda)\Q_2\in \mathcal{G}_1$ for any $\lambda\in[0,1]$.\\ 
Since $D_1,\ D_2$ are both greater than zero, $D_{\lambda}:=\lambda D_1+(1-\lambda) D_2>0$. It similarly follows that  $W_{\lambda}:=\lambda W_1+(1-\lambda) W_2\geq0$. Let us now verify that $W_{\lambda}/D_{\lambda}\leq 1$, this is equivalent to verify that:
\begin{equation}
	\label{convexityW}
	\lambda (W_1-D_1)+(1-\lambda) (W_2-D_2)\leq 0,
\end{equation}
since $D_{\lambda}>0$.
By hypothesis $W_i/D_i\leq1$ this means that $W_i-D_i \leq 0$ for $i\in\{1,2\}$. Thus equation (\ref{convexityW}) is immediately satisfied.
We now only need to show that $g(\Q_{\lambda})>0$:
\begin{align}
	E_{\lambda}:&=\lambda E_1 +(1-\lambda)E_2>\lambda \sqrt{(D_1+qW_1)^2+S_1^2}+ (1-\lambda) \sqrt{(D_2+qW_2)^2+S_2^2}\nonumber\\
	&\geq \sqrt{(\lambda(D_1+qW_1)+(1-\lambda)(D_2+qW_2))^2 + (\lambda S_1 + (1-\lambda)S_2)^2}\nonumber\\
	&\mathrel{=} \sqrt{(D_{\lambda}+qW_{\lambda} )^2+S_{\lambda}^2}\,,
\end{align}
where the Minkowski inequality was used in the inequality computation.\\
Let us now consider the auxiliary function 
\begin{equation}
	\label{eq:barPhi}
	\Bar{\Phi}(p)=(	\Gamma -1)^2 (h_1(p)-h_2(p))
\end{equation}
where
\begin{align}
	h_1(p)&=\left(S^2+(E+p)\left(\frac{p}{\Gamma-1}-E\right)\right)^2\\
	h_2(p)&=(D+qW)^2((E+p)^2-S^2).
\end{align}
The expressions $h_1(p)$ and $h_2(p)$ above are obtained by multiplying
	Eq.~\eqref{eq:zerophip} by $E+p$ and then squaring both sides.
	Consequently, every zero of $\Phi(p)$ is also a zero of
	$\bar{\Phi}(p)$, while the converse does not hold in general.
 The function $\Bar{\Phi}(p)$ can be written as:
\begin{equation}
	\label{eq:barPhi4pol}
	\Bar{\Phi}(p)=c_0+c_1p+c_2p^2+c_3p^3+c_4p^4
\end{equation}
where
\begin{align}
	c_0&=(\Gamma-1)^2(S^2-E^2)(S^2-E^2+(D+qW)^2)\\
	c_1&=2E(\Gamma-1)(2-\Gamma)(S^2-E^2)-2E(D+qW)^2(\Gamma-1)^2\\
	c_2&=E^2(6-6\Gamma+\Gamma^2)+2S^2(\Gamma-1)-(D+qW)^2(\Gamma-1)^2\\
	c_3&=2E(2-\Gamma)\\
	c_4&=1.
\end{align}
The nature of the roots of  
the polynomial $\Bar{\Phi}(p)$ has already been analyzed by \cite{Cai2024Newton}, and their reasoning 
expressed by Lemmas (2.1), (2.2), and (2.3)
remains valid also in the presence of reactions. In particular, it is possible to show that there are either
two positive and two negative roots or two positive and two complex roots.
Similarly,  the smallest positive root of $\Bar{\Phi}(p)$ is the unique positive root of $\Phi(p)$, which is the physical pressure $p(\Q)$.
A first possible approach is to find $p^*$ as the solution of $\Bar{\Phi}(p^*)=0$ using a Newton algorithm and initial guess as in Algorithm 2.1 of \cite{Cai2024Newton}.
Alternatively, we adopt:
\begin{equation}
	\psi(p):=(E+p)\Phi(p)=S^2+(E+p)\left(\frac{p}{\Gamma-1}-E\right)+ (D+qW)\sqrt{(E+p)^2-S^2}.
\end{equation}
with an initial guess $p_c^\Q$ given by
\begin{equation}
	p_c^\Q:=\frac{1}{2}\left( (\Gamma-2)E+ \sqrt{(2-\Gamma)^2E^2-4(\Gamma-1)\left((S^2-E^2)+\left(D+qW\right)\sqrt{E^2-S^2}\right)}\right)
\end{equation}
such that $h_1(p_c^\Q)=h_2(0)$ (see Lemma (2.6) of \cite{Cai2024Newton}).
Hence, on a practical ground we recover the unknown pressure 
$p^*$  as the solution of $\psi(p)=0$ using a Newton algorithm with an initial guess $p^0=0$ if $(D+qW)\geq \frac{E^2-S^2}{E}$ or   $p^0=p_c^\Q$ otherwise. Once the pressure has been computed, the velocity and density follow immediately from Eq.~\eqref{eq:rec-v}-\eqref{eq:rec-rho}.

\subsection{Treatment of the reaction source term}
The combustion process is characterized by fast reacting components that reach their final equilibrium almost instantly. The source term needs to take into account the instant reaction, hence it is considered to be stiff. This rapid change in the reactant creates numerical difficulties in the computation of the solution of the evolutionary problem.\\
To avoid instability, it is common to use implicit schemes in case of stiff problems.
In our numerical simulations we employ a classical splitting technique, hence we first assume all source terms to be zero. In these conditions we compute a predictor $\Q^{*}$ explicitly from the homogeneous hyperbolic system of conservation laws (\ref{eq:generalALEscheme}). We then use this predictor solution as the initial guess of a Newton iteration that accounts only for the reaction source. Namely, we look for  $\tilde{\Q}$ such that 
\begin{equation}
	\tilde{\Q}-\Q^*-\Delta t \S(\tilde{\Q})=0\,.
\end{equation}
From our numerical experiments we find that in general around 3 iterations are needed for the convergence of the Newton algorithm. Although the homogeneous ALE discretization and its TVD reconstruction are second-order accurate, the Lie splitting and the backward-Euler treatment of the reaction source make the complete reactive balance-law solver first-order accurate in time in general.

\section{Numerical tests}
\label{sec.tests}
In this section we present some numerical tests that illustrate the propagation of different relativistic detonation waves.
The numerical solutions are obtained with the ALE finite volume scheme 
described above,
using a moving mesh with a constant velocity equal to the shock speed $V_s$. The numerical solution obtained with the ALE finite volume scheme is then validated against the reference solution obtained from the corresponding ODE, described in Sect.~\ref{sec:ZNDprofile} and indicated in the following as the \emph{exact} solution as we can solve the ODE up to arbitrary precision.

\noindent In all of our tests we consider the spatial domain $\Omega=[-1,1]$ where at  $\xi=0$ lies the discontinuity, and Dirichlet boundary conditions were considered.
In test cases \textit{1, 2, 3} we have considered $N = 10880$ as the number of cells and  CFL factor given by $CFL=0.5$.  For ease of computation we have employed in all of our numerical tests the Rusanov flux, since the Osher-type flux \cite{OsherUniversal} requires the evaluation of the eigenvector matrix, which in the case of relativistic flows is purely numerical and therefore computationally more demanding.
As a general criterion, we state that the ALE scheme converges to a stationary solution in the comoving frame if $\lVert\Q^{n+1}-\Q^{n}\rVert_2/\lVert\Q^{n}\rVert_2<\epsilon$ for a set tolerance  $\epsilon\sim 10^{-12}$. An important remark for a correct comparison among the ODE solution and the ALE one is the following: the ODE is actually solved in the comoving frame of the shock through an explicit Lorentz transformation, cf.~\eqref{eq:lorentz}. The ALE scheme, on the contrary, follows the shock front on a moving mesh but with respect to the rest frame, hence it does not imply any Lorentz contraction. Thus, the ZND profile has  maximum length in the comoving frame, while, for a proper comparison, the width of all  ALE profiles must be expanded  through multiplication by $\gamma_s$.

\subsection{Numerical convergence test}
\begin{figure}[ht]
	\centering
	
	\begin{minipage}[t]{0.47\textwidth}
		\vspace{0pt}
		\centering
		\includegraphics[width=\linewidth]{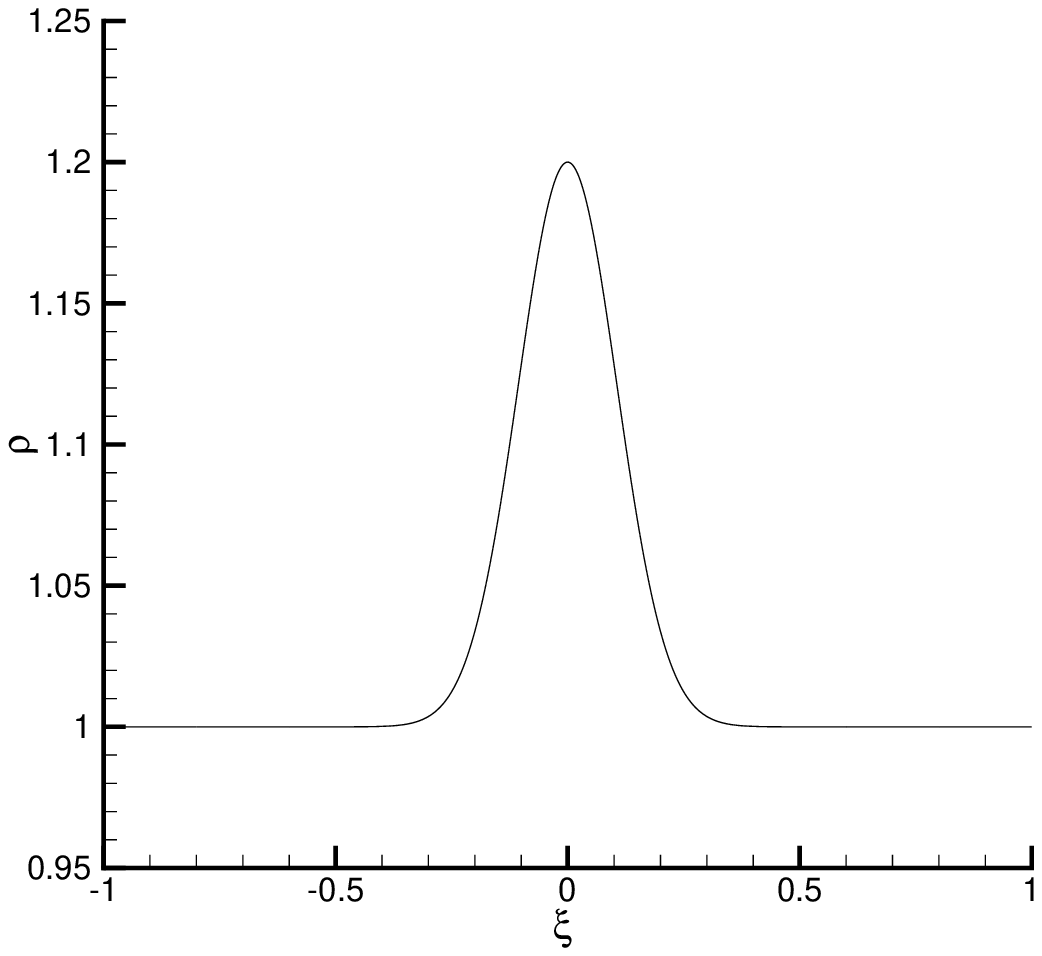}
		\caption{Initial smooth density profile used in the convergence test.}
		\label{fig:smooth_density_profile}
	\end{minipage}
	\hfill
	\begin{minipage}[t]{0.49\textwidth}
		\vspace{0pt}
		\centering
		\small
		\setlength{\tabcolsep}{4pt}
		\renewcommand{\arraystretch}{1.05}
		\begin{tabular}{lll}
			\toprule
			\textbf{Category}
			& \textbf{Quantity}
			& \textbf{Value} \\
			\midrule
			Physical parameters
			& $\Gamma$   & $4/3$   \\
			& $K_0$      & $20.0$  \\
			& $q$        & $0.8$   \\
			& $T_i$      & $2.55$  \\
			& $m$        & $0.0$   \\
			& $\gamma_s$ & $2.931$ \\
			& $E_a$      & $2.0$   \\
			\midrule
			Smooth profile
			& $\rho_0$   & $1.0$   \\
			& $p_0$      & $1.0$   \\
			& $v_0$      & $0.94$  \\
			& $Z_0$      & $0$     \\
			& $A$        & $0.2$   \\
			& $\sigma$   & $0.15$  \\
			\bottomrule
		\end{tabular}
		\captionof{table}{Physical parameters and initial data of the smooth-density test.}
		\label{tab:smooth_density_test}
	\end{minipage}
\end{figure}

In order to verify the order of convergence of our numerical scheme of the homogeneous problem, we choose a constant pressure, i.e. $p(x,t)=p_0$, constant speed, i.e. $v(x,t)=V$ and smooth density profile, i.e. $\rho(x,t)=\phi(\xi)$, where $\xi=x-Vt$ and
$\phi$ a generic smooth positive function.  As such, the resulting profile propagates un-distorted at the constant speed $V$ and it obeys the relativistic Euler equations, as it is easy to check.
As a possible density profile we choose
\begin{equation}
	\rho(x,t)=\rho_0\left[1+A
	\exp\left(-((x-Vt)/\sigma)^2\right)\right],
	\qquad |A|<1.
\end{equation}
In our numerical tests we use as physical variables the one represented in Table~\ref{tab:smooth_density_test}. Due to the particular choice of the ignition temperature $T_i$, in this test case the source term vanishes identically, so that no operator-splitting or implicit source-integration error is introduced; therefore, second-order convergence of the homogeneous ALE scheme is expected. However we formally still solve the complete reactive system, as in all subsequent tests. 
For the convergence analysis, all simulations are stopped at the same prescribed physical time, $t=1$, and the numerical solution is compared with the exact translated profile. In particular, we do not use the stationary-state stopping criterion $\lVert\Q^{n+1}-\Q^{n}\rVert_2/\lVert\Q^{n}\rVert_2<\epsilon$
in this test. Although the exact solution is stationary in the comoving frame, applying this criterion separately at each resolution would generally result in different final times and could make the measured errors depend on the prescribed tolerance. This may produce artificially large or otherwise unreliable convergence rates once the discretization error becomes comparable to the stopping tolerance or to round-off errors.
Table~\ref{tab:smooth_density_convergence} reports the errors and observed convergence rates for the rest-mass density. 
The $L_1$ error displays an approximately second-order convergence rate, consistently with the second-order space-time reconstruction used in the homogeneous ALE scheme. The velocity, pressure, and Lorentz factor remain constant up to machine precision throughout the computation. The mass fraction of the reactant $Z$ is also preserved exactly as a constant. The temperature and the conserved variables that depend on the density exhibit convergence behaviour consistent with that of $\rho$.
The initial smooth density profile used in this convergence test is shown
in Fig.~\ref{fig:smooth_density_profile}.
\begin{table}[!tb]
	\centering
	\caption{$L_1$ errors and observed convergence rates for the
		rest-mass density $\rho$ at the fixed final time $t=1$.}
	\label{tab:smooth_density_convergence}
	\begin{tabular}{rcc}
		\toprule
		$N$ & $L_1$ error & Order \\
		\midrule
		400 & $1.0776\cdot 10^{-4}$ & --   \\
		800 & $2.9533\cdot 10^{-5}$ & 1.87 \\
		1600 & $7.9502\cdot 10^{-6}$ & 1.89 \\
		3200 & $2.0983\cdot 10^{-6}$ & 1.92 \\
		6400 & $5.4661\cdot 10^{-7}$ & 1.94 \\
		\bottomrule
	\end{tabular}
\end{table}
We also tested the robustness of the scheme under extreme low-density and low-pressure conditions. We test our numerical scheme on the above homogeneous problem where all parameters  are fixed as in Table~\ref{tab:smooth_density_test}.
	To examine performance at low densities, we decrease  $\rho_0$ for the low-density test. Stable computations were obtained for rest-mass densities as small as $\rho_0=10^{-10}$. For the low-pressure test, instead we study the results by decreasing the value of $p_0$. The smallest pressure successfully tested is $p_0=10^{-5}$, whereas $p_0=10^{-6}$ led to the generation of a nonphysical state. At very low pressure the recovery of the pressure is increasingly sensitive to round-off errors and loss of significance. These limit values represent the robustness limit for the present test configuration, double-precision arithmetic, and selected numerical tolerances. 
\subsection{Test 1}
\begin{figure}[!htbp]
	\begin{center}
		\begin{tabular}{cc}
			\includegraphics[width=0.45\textwidth]{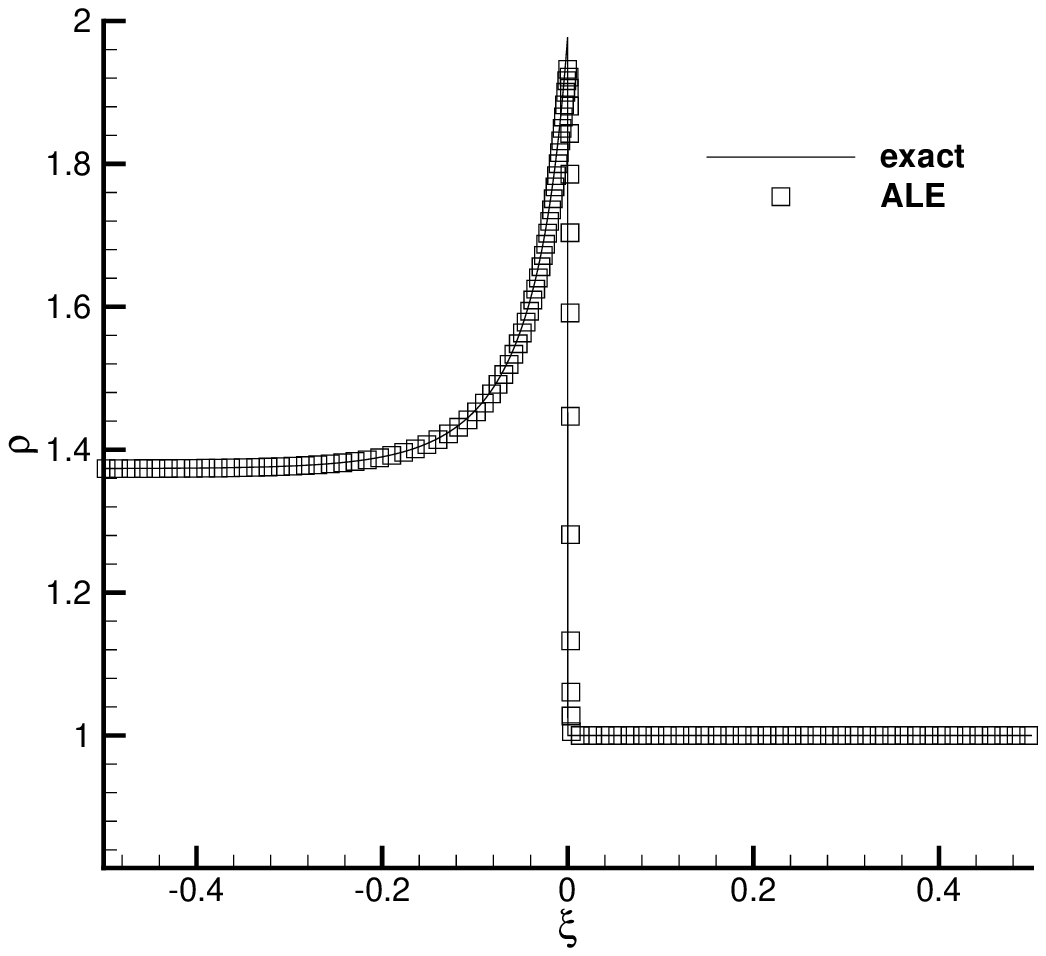} &
			\includegraphics[width=0.45\textwidth]{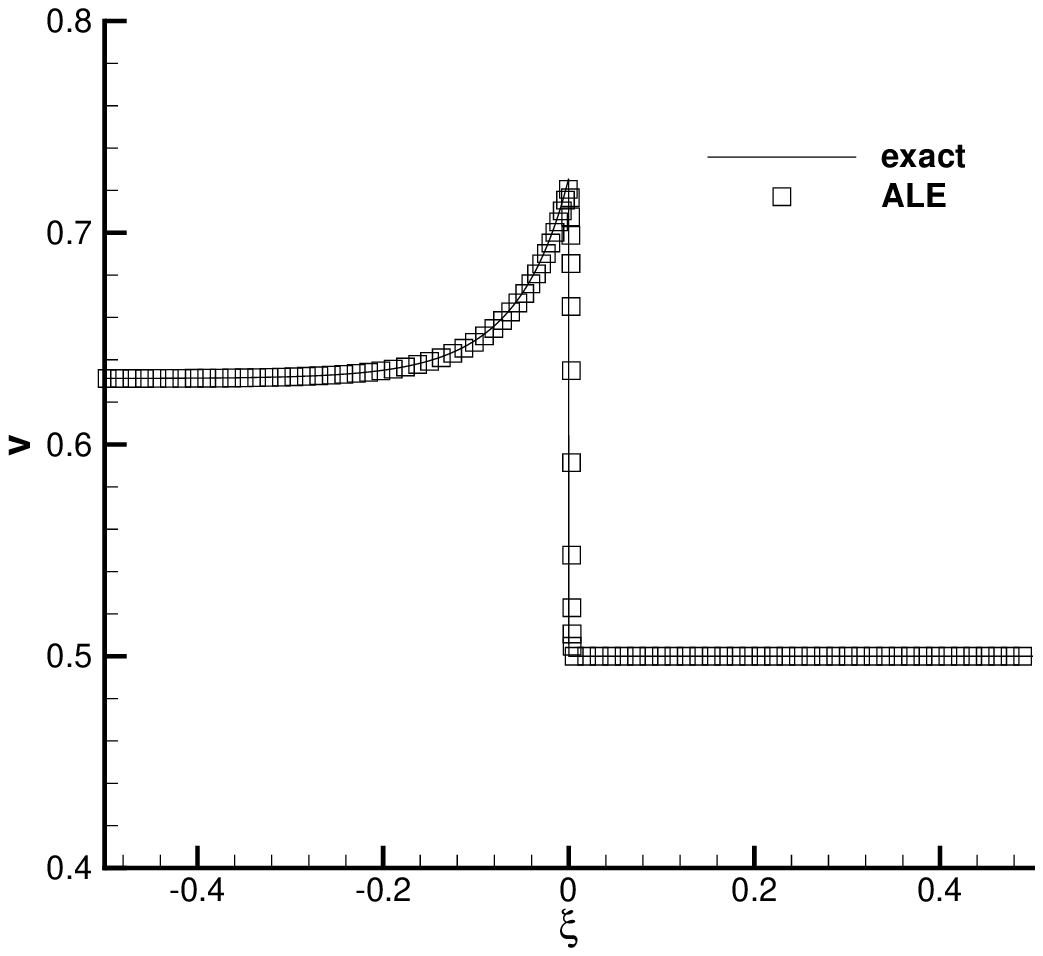} \\
			\includegraphics[width=0.45\textwidth]{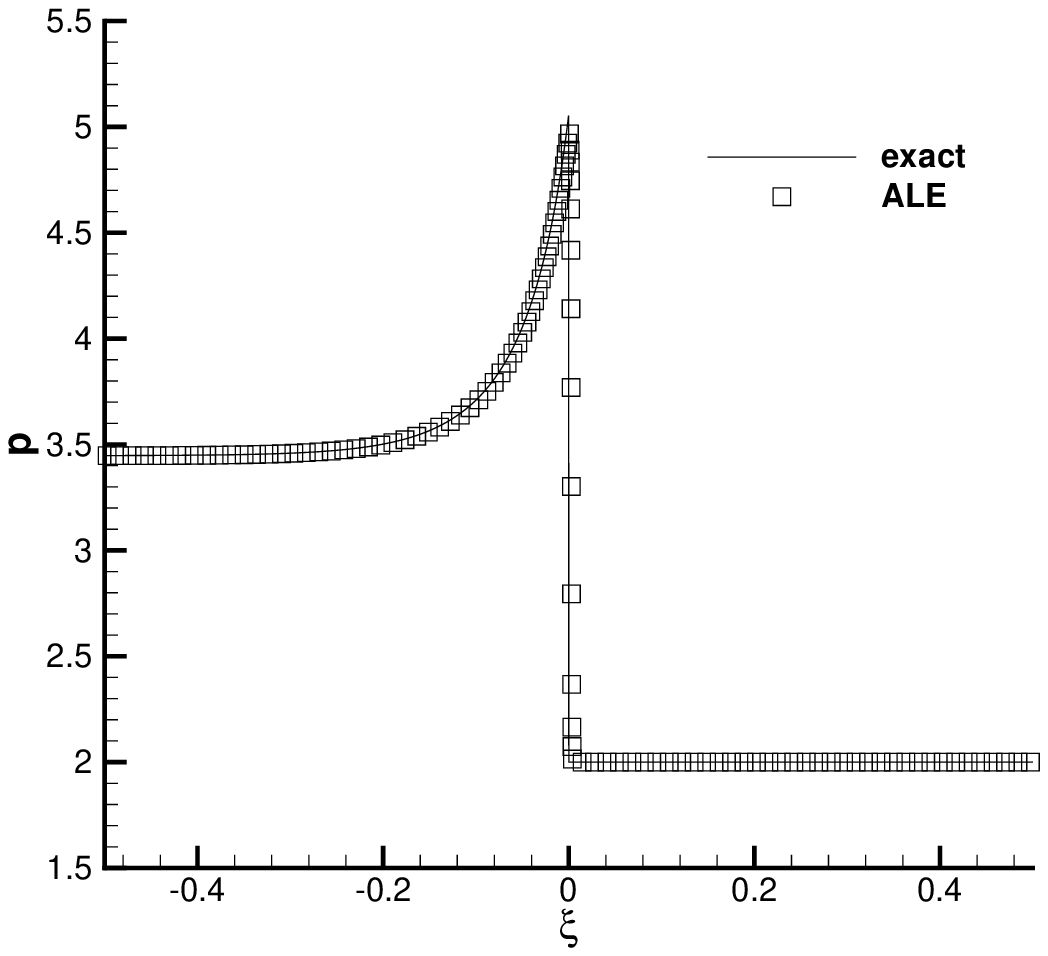} &
			\includegraphics[width=0.45\textwidth]{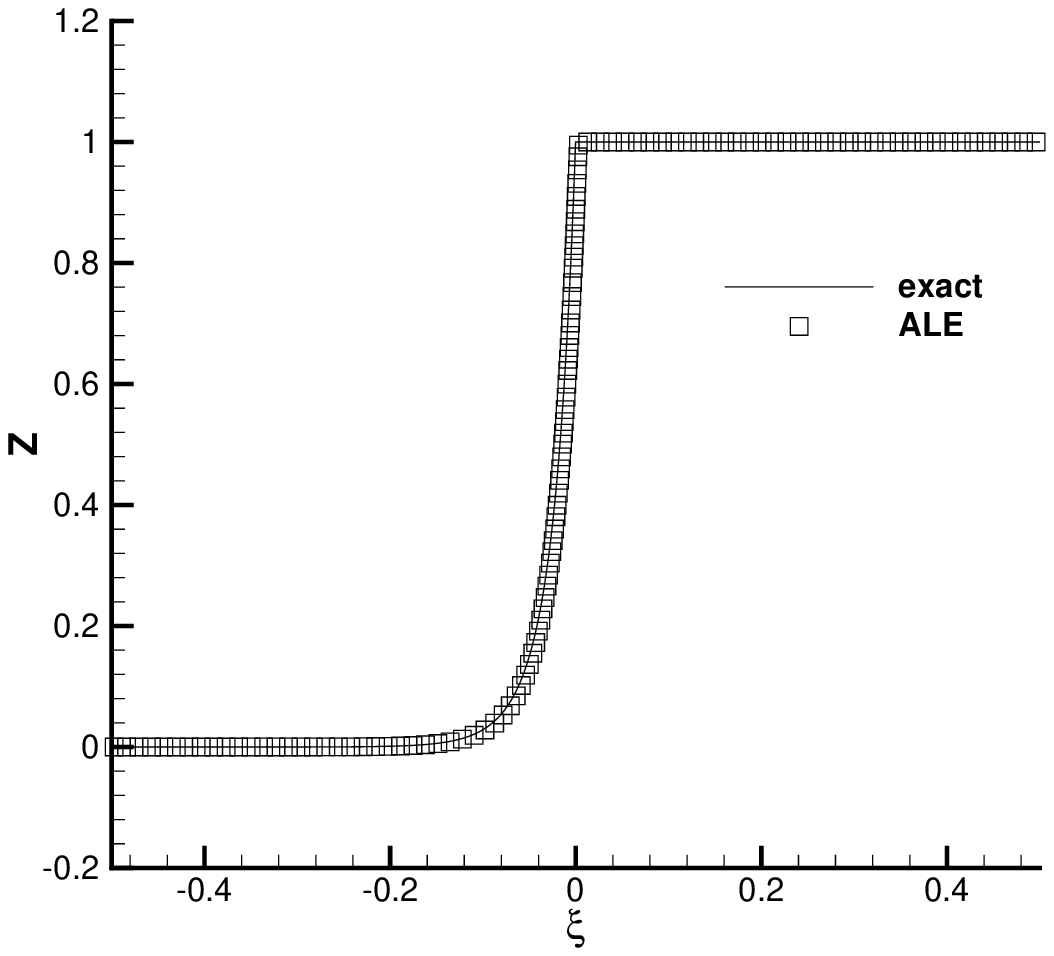} \\
			\includegraphics[width=0.45\textwidth]{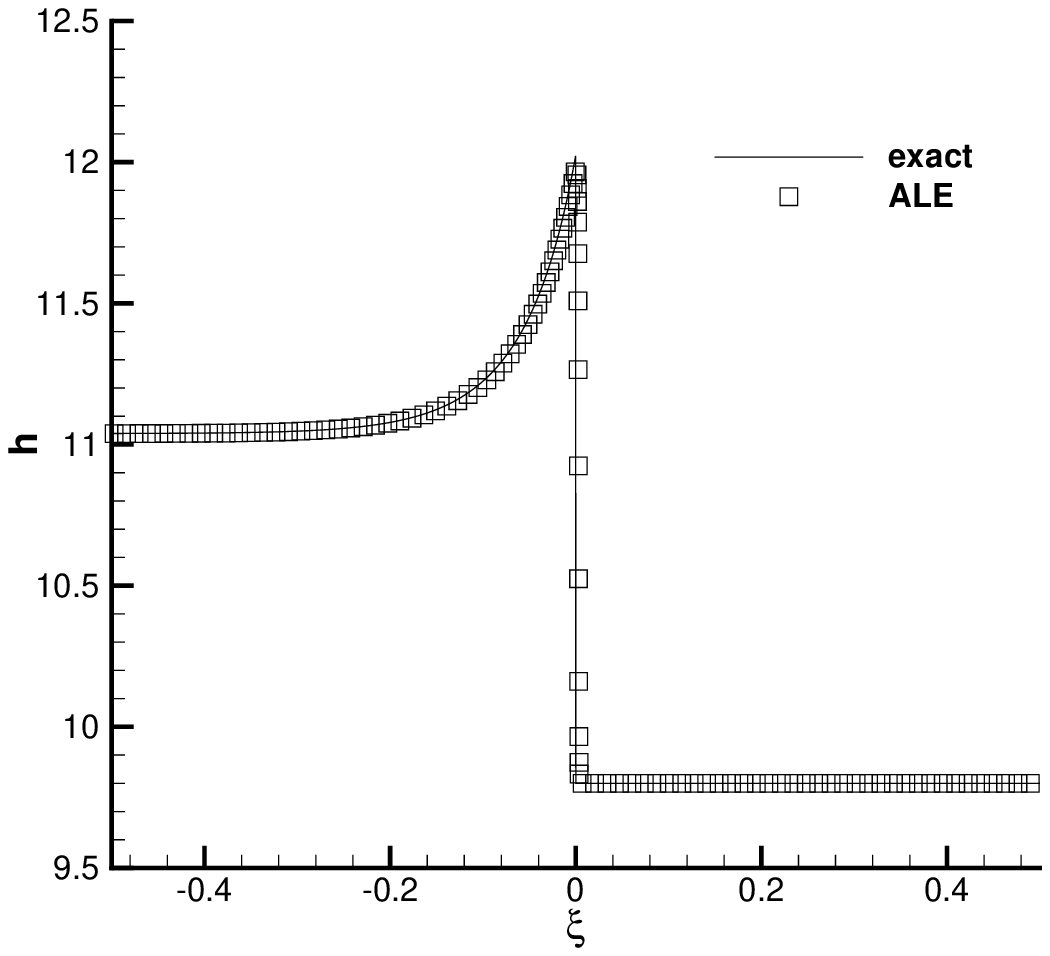} &
			\includegraphics[width=0.45\textwidth]{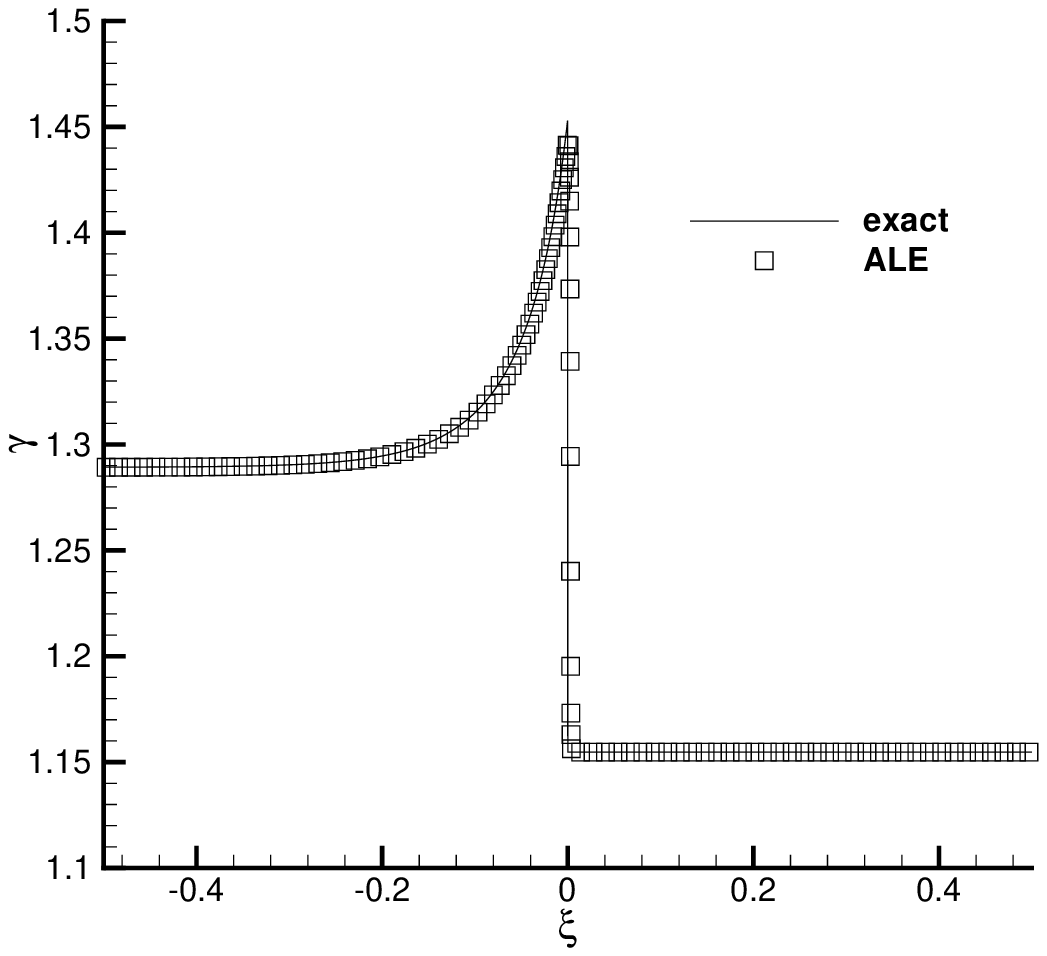} \\
		\end{tabular}
		\caption{(Test 1): Numerical solution of the relativistic detonation for a CJ mass flux obtained with the ALE scheme, compared to the quasi exact solution given by the ODE solver.}
		\label{fig.test1}
	\end{center}
\end{figure}

The physical setting in this test case is given by the quantities in Tab.~\ref{table1}. This test simulates the shock profile for the case of a Chapman-Jouguet mass flux.
\begin{table}[!tb]
	\centering
	\caption{Physical parameters and initial data of Test 1.}
	\label{table1}
	\begin{tabular}{lll}
		\toprule
		\textbf{Category} & \textbf{Quantity} & \textbf{Value} \\
		\midrule
		Physical parameters & $\Gamma$ & $4/3$ \\
		& $K_0$    & 20.0 \\
		& $q$    & 0.8 \\
		& $T_i$ & 2.25 \\
		& $m_{CJ}$      & 0.9059044176708 \\
		& $\gamma_s$    & 2.0810839845077 \\
		\midrule
		Unshocked state     & $\rho_0$ & 1.0 \\
		& $p_0$    & 2.0 \\
		& $v_0$    & 0.5 \\
		& $Z_0$    & 1.0 \\
		\midrule
		Shocked burnt state    & $\rho_b$ & 1.3735273312153 \\
		& $p_b$    & 3.4470794604194 \\
		& $v_b$    & 0.6311617316827 \\
		& $Z_b$    & 0.0 \\
		\bottomrule
	\end{tabular}
\end{table}
The source term is modeled with the discrete ignition temperature model of Eq.~\eqref{Eq:K-ignition}. In this configuration the difference between the shocked burnt state and the shocked unburnt state is relatively small. 
The results are shown in Fig.~\ref{fig.test1} and they prove the ability of the ALE scheme in reproducing the 
Zel'dovich spike with a very good accuracy. We also stress that the ALE algorithm uses as initial conditions the unshocked state and the shocked fully burnt state, while the shocked unburnt state is recovered as a result of the simulation.
The choice of $T_i$ is critical. It has to satisfy the constraint $T_0<T_i<T_b$. In this test case $T_0=2.0$ and $T_b=2.509655$. Numerical experiments indicate that selecting the ignition temperature $T_i$ closer to the initial temperature $T_0$ results in smaller spike amplitudes, while choosing $T_i$ too close to the shocked burnt state temperature leads to increased propagation of numerical viscosity which, ultimately, leads to incorrect solutions.

\subsection{Test 2}
\begin{figure}[!htbp]
	\begin{center}
		\begin{tabular}{cc}
			\includegraphics[width=0.45\textwidth]{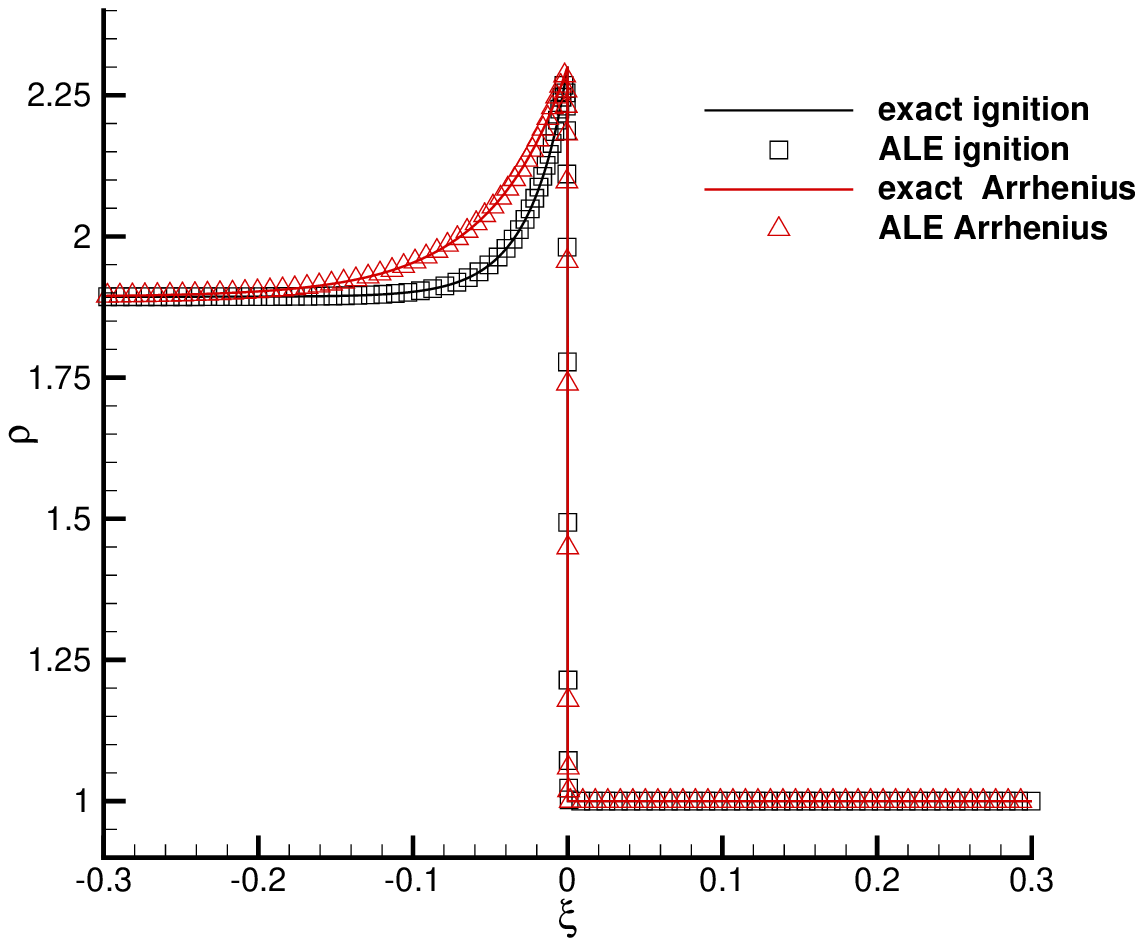} &
			\includegraphics[width=0.45\textwidth]{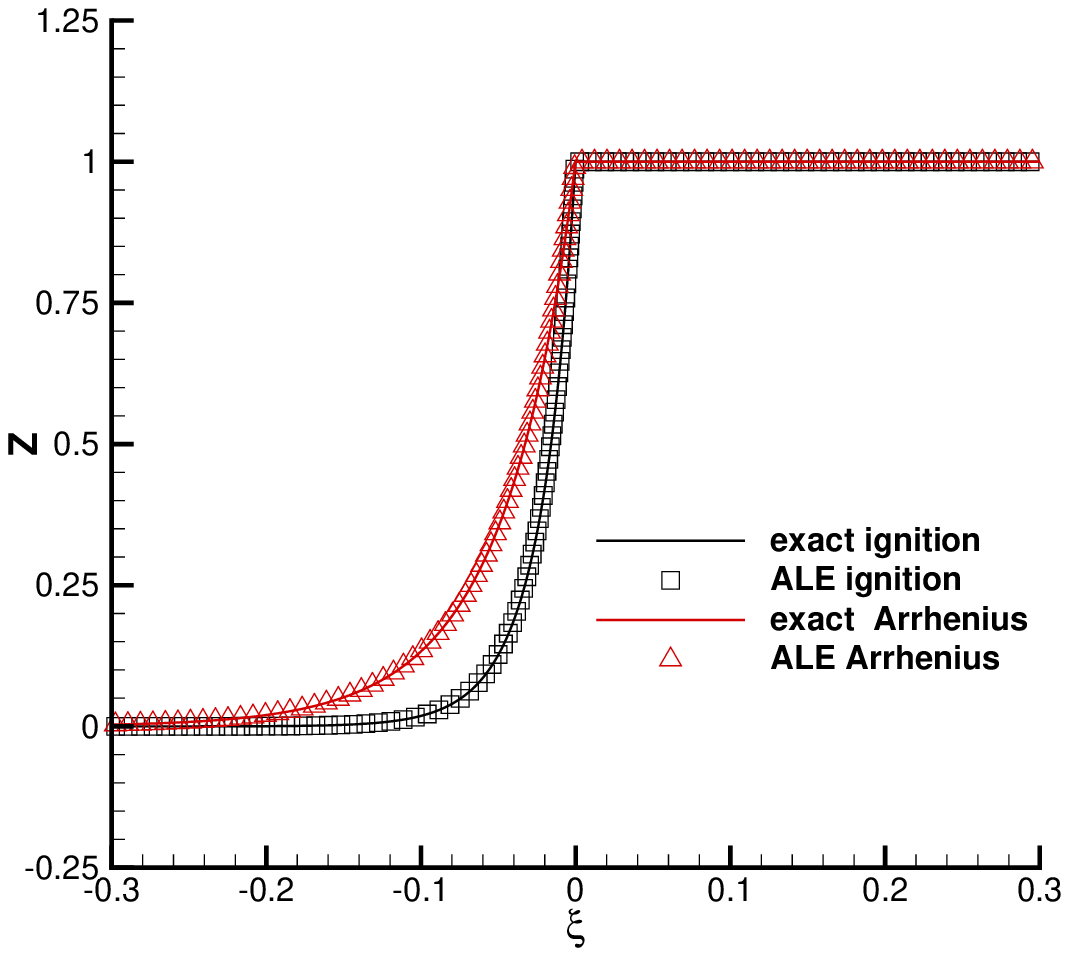} \\
		\end{tabular}
		\caption{(Test 2): Numerical solutions of the relativistic detonation waves obtained with the ALE scheme, compared to the quasi exact solution given by the ODE solver. The different profiles are obtained with the discrete ignition model for the source term (black line) and with the Arrhenius law (red line).}
		\label{fig.test2}
	\end{center}
\end{figure}

The aim of this second test is to investigate the effect of the source term on the ZND profile, in particular 
the difference between a sharp activation of the combustion, see Eq.~\eqref{Eq:K-ignition}, and the one provided by Arrhenius law, see Eq.~\eqref{Eq:K-Arrhenius}.
The physical setting in this test case is given by the quantities in Table~\ref{table2}.  The corresponding results are shown in Fig.~\ref{fig.test2} only for the rest mass density and the mass fraction of the chemical reactant. For the numbers chosen, the combustion process is faster with abrupt ignition. For smaller values of the activation energy $E_a$, though, the combustion can be made increasingly faster. Similarly, sharper profiles are obtained for larger values of $K_0$.
\begin{table}[!tb]
	\centering
	\caption{Physical parameters and initial data of Test 2.}
	\label{table2}
	\begin{tabular}{lll}
		\toprule
		\textbf{Category} & \textbf{Quantity} & \textbf{Value} \\
		\midrule
		Physical parameters & $\Gamma$ & $4/3$ \\
		& $K_0$    & 20.0 \\
		& $q$    & 0.8 \\
		& $T_i$ & 2.55 \\
		& $m$      & 1.0 \\
		& $\gamma_s$ &1.4142135623731\\
		& $E_a$    & 2.0 \\
		\midrule
		Unshocked state     & $\rho_0$ & 1.0 \\
		& $p_0$    & 2.0 \\
		& $v_0$    & 0.0 \\
		& $Z_0$    & 1.0 \\
		\midrule
		Shocked burnt state    & $\rho_b$ &  1.8932771241913 \\
		& $p_b$    & 5.3271556030358\\
		& $v_b$    & 0.3584408321366\\
		& $Z_b$    & 0.0 \\
		\bottomrule
	\end{tabular}
\end{table}

\subsection{Test 3}
\begin{figure}[!htbp]
	\begin{center}
		\begin{tabular}{cc}
			\includegraphics[width=0.45\textwidth]{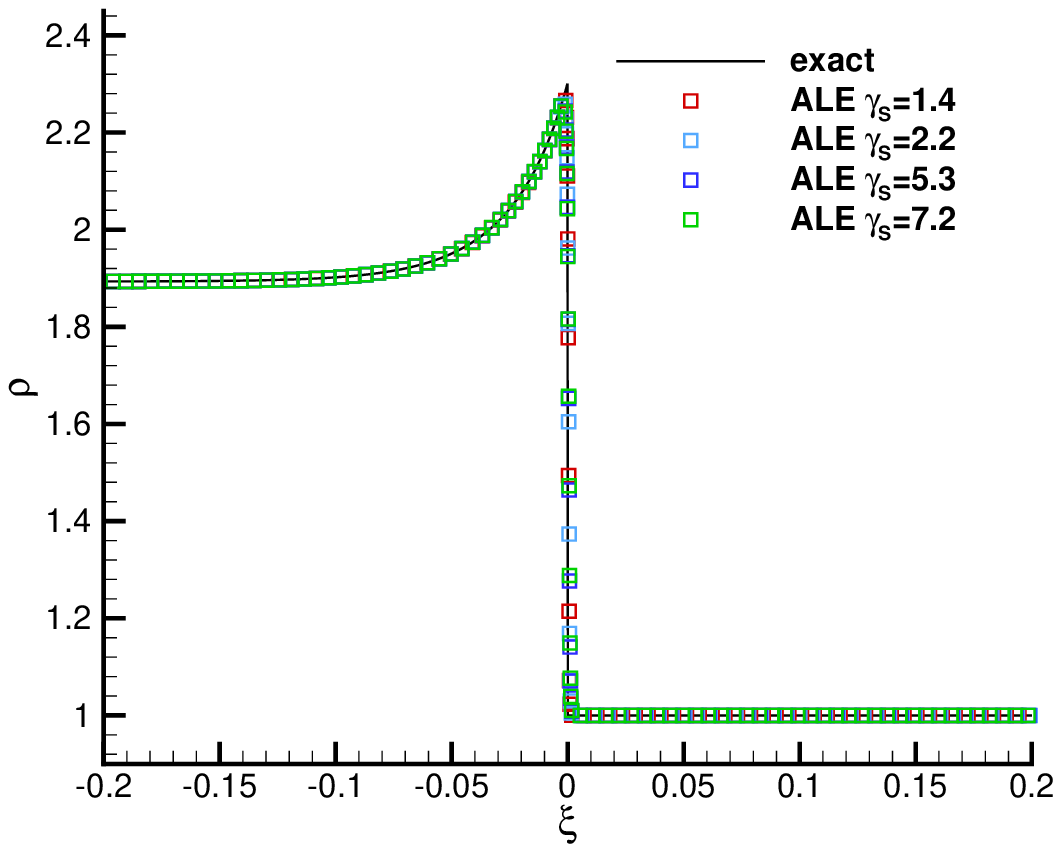} &
			\includegraphics[width=0.45\textwidth]{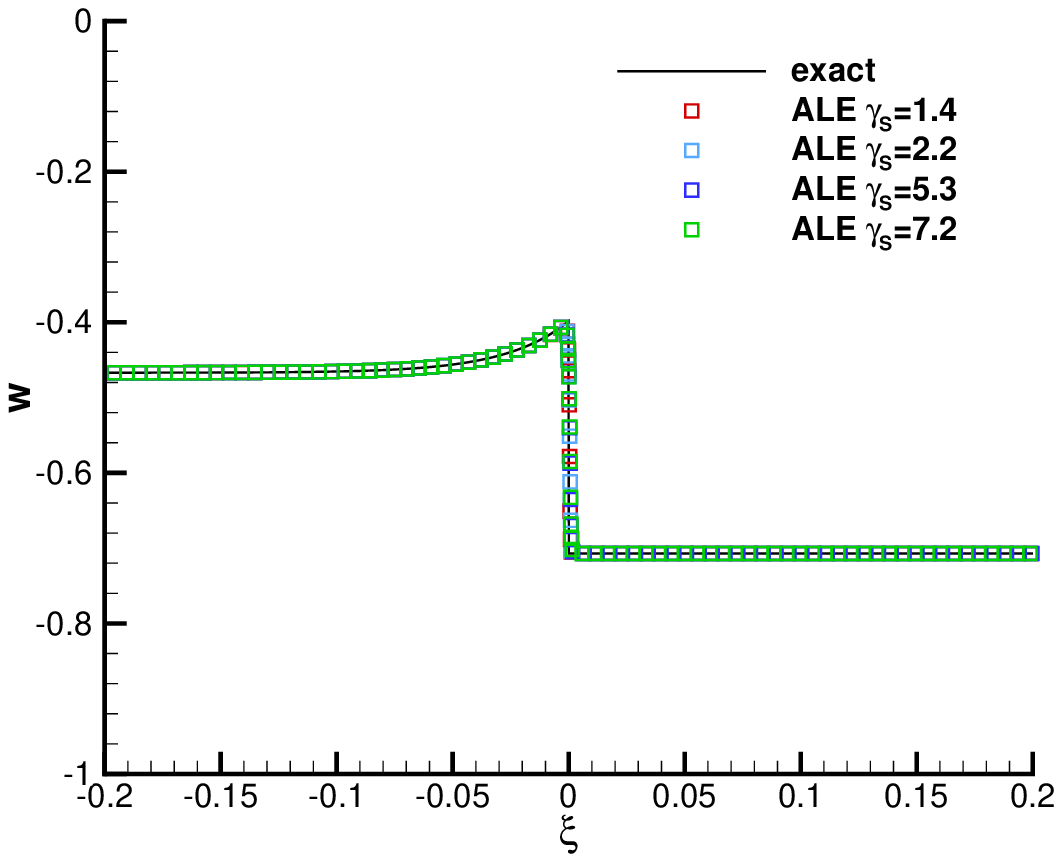} \\
			\includegraphics[width=0.45\textwidth]{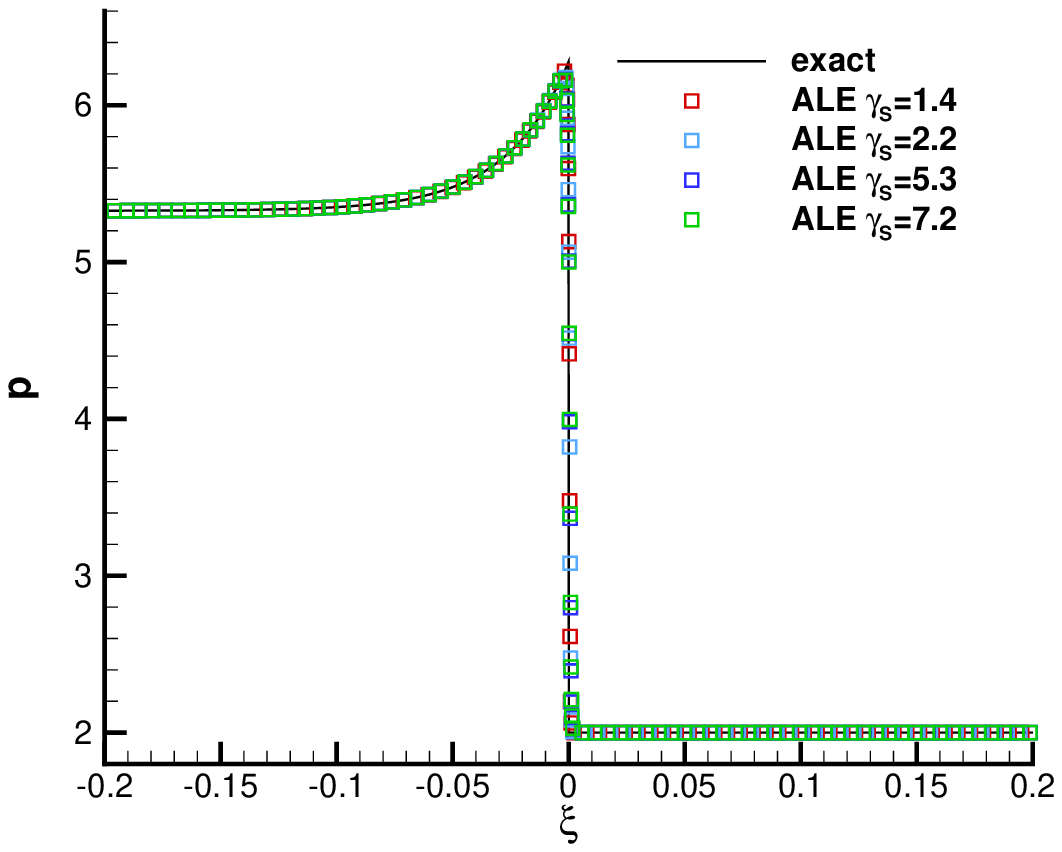} &
			\includegraphics[width=0.45\textwidth]{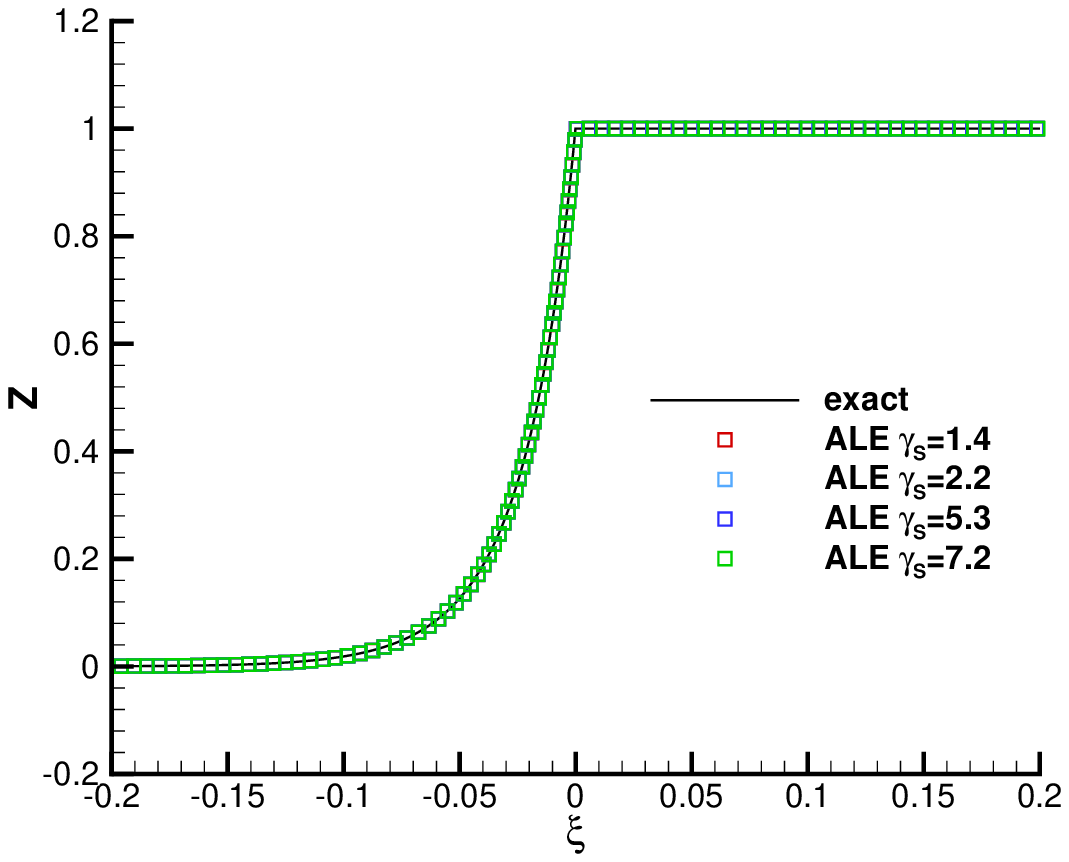} \\
			\includegraphics[width=0.45\textwidth]{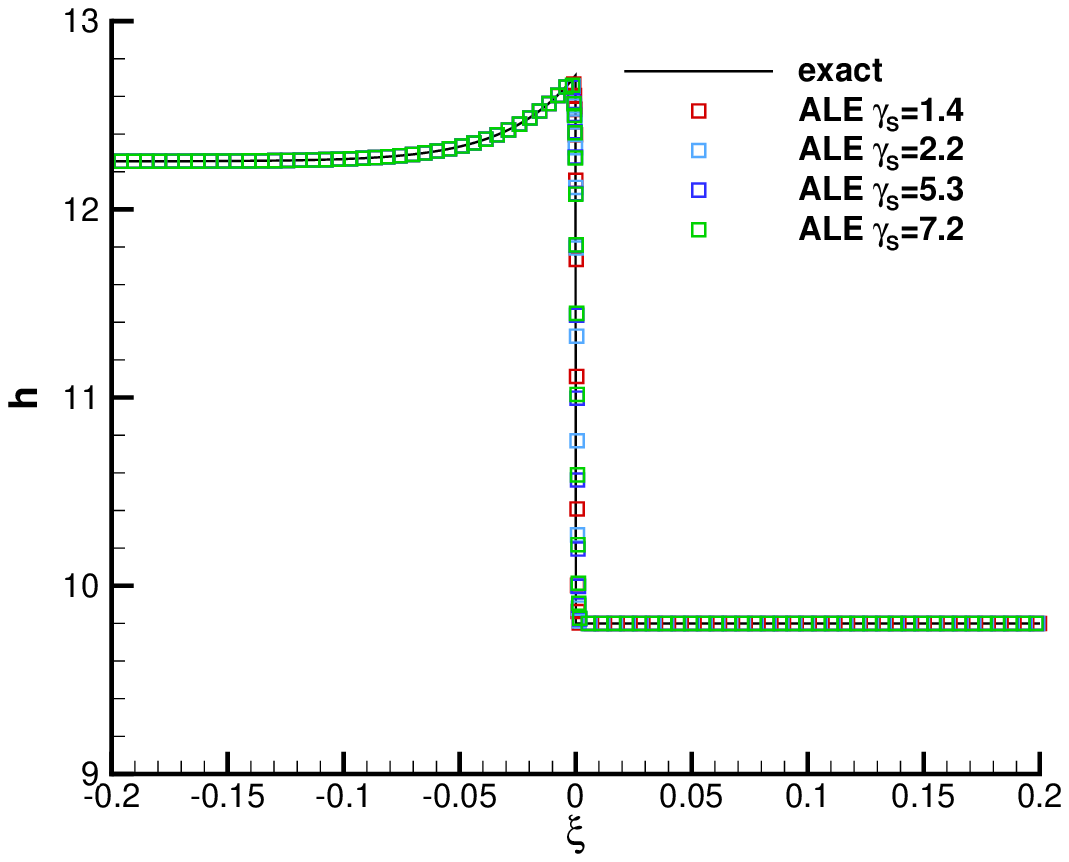} &
			\includegraphics[width=0.45\textwidth]{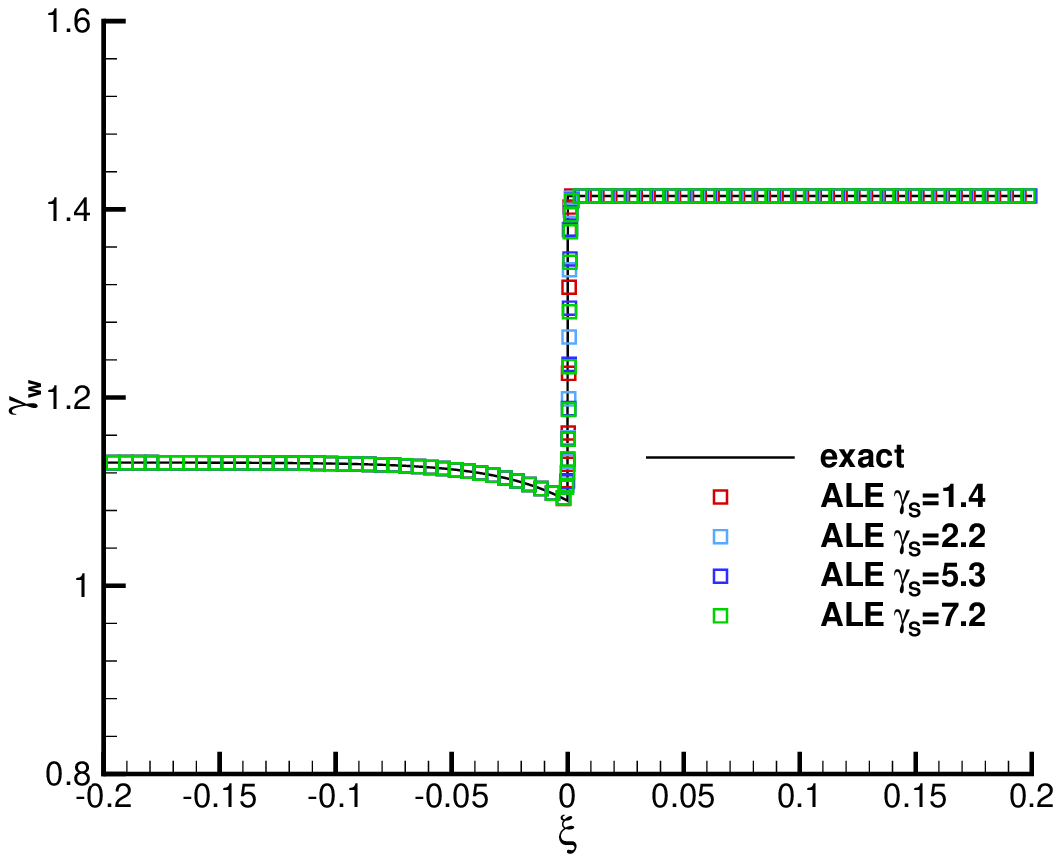} \\
		\end{tabular}
		\caption{(Test 3): Numerical solution of the relativistic detonation of Test 3 when increasing values of the Lorentz factor $\gamma_s$ are considered.}
		\label{fig.test3}
	\end{center}
\end{figure}

With this test, we check the ability of the code in dealing with larger and larger velocities of the shock front.
The physical settings are given in Table~\ref{table3}. In this case we set the mass flux $m=1.0$ and we are interested to evaluate the solution for different values of $\gamma_s$. To this extent it is necessary to increase $v_0$. This follows directly from the definition of the mass flux, which can be computed as $m=-\gamma_{w,0}\rho_0 w_0$.
By keeping $m$ and $\rho_0$ equal for all tests,  $\gamma_{w,0}$ and $w_0$ remain unchanged. 
Hence, the Lorentz factor
$\gamma_s$ is fully controlled by the value of $v_0$.
Table~\ref{table3} contains all this information for different values of the initial velocity.
After correcting for the Lorentz contraction, we expect all the profiles to match each other. This
is indeed the case as reported in Fig.~\ref{fig.test3} where the exact solution is again reported for comparison.
We stress that, in order to cope with large velocities, it has been fundamental to use the algorithm shown in Sect.~\ref{sec:cons2prim} for the conversion from the conservative to the primitive variables.

\begin{table}[!tb]
	\centering
	\caption{Physical parameters and initial data of Test 3.}
	\label{table3}
	\begin{tabular}{lll}
		\toprule
		\textbf{Category} & \textbf{Quantity} & \textbf{Value} \\
		\midrule
		Physical parameters & $\Gamma$ & $4/3$ \\
		& $K_0$    & 20.0 \\
		& $q$    & 0.8 \\
		& $T_i$ & 2.55 \\
		& $m$      & 1.0 \\
		\midrule
		Unshocked state     & $\rho_0$ & 1.0 \\
		& $p_0$    & 2.0 \\
		& $Z_0$    & 1.0 \\
		\midrule
		Shocked burnt state    & $\rho_b$ &  $1.8932771241913$\\
		& $p_b$    & $5.3271556030358$ \\
		& $Z_b$    & 0.0 \\
		\midrule 
		Case 1    & $\gamma_s$ &  $1.4142135623731$ \\
		          & $v_0$      &  $0.0$   \\    
		          & $v_b$      &  $0.3584408321366$   \\    
		\midrule 
		Case 2    & $\gamma_s$ &  $2.2103434310451$ \\
		& $v_0$      &  $0.5$    \\
        & $v_b$      &  $0.7279731765489$   \\    
		\midrule 
		Case 3    & $\gamma_s$ &  $5.3091700274503$ \\
		& $v_0$      &  $0.9$     \\
        & $v_b$      &  $0.9514924584242$   \\    
		\midrule 
		Case 4    & $\gamma_s$ &  $7.2131690172665$ \\
		& $v_0$      &  $0.945$    \\
        & $v_b$      &  $0.9736422996363$   \\    
		\bottomrule
	\end{tabular}
\end{table}

\clearpage
\section{Conclusions}
\label{sec.conclusions}

We have used a one-dimensional 
Arbitrary-Lagrangian-Eulerian (ALE) finite volume scheme to solve the dynamics of relativistic detonation waves. Apart from an additional equation for the
mass fraction of the reactant, the remaining equations preserve the same formal
structure as those of an inert shock, since the chemistry of the process is accounted for
by a simple extra term in the definition of the total energy density. Our ALE scheme is based on a combination of the following main properties:
\begin{itemize}
	\item flux correction to track the propagation of the shock front in its rest frame;
	\item TVD reconstruction in space;
	\item usage of a particularly robust algorithm for the recovering of the primitive variables from the conserved ones, adapted from \cite{Cai2024Newton};
	\item operator splitting to treat the stiff reaction source terms in a robust manner. 
\end{itemize}
In this way we have been able to reproduce the relativistic Zel'dovich-von Neumann-Doering profile with very good  accuracy, up to Lorentz factors of the shock front $\gamma_s\sim 7$.

In addition, we have provided a clean strategy for the numerical computation of the Chapman-Jouguet mass flux, which was so far missing in the relativistic regime.
When compared to the Newtonian regime, relativistic detonations manifest a sharp qualitative difference. In the Newtonian case, the Zel'dovich pressure jump $\Delta p_Z$, i.e. the pressure difference between the shocked unburnt and fully burnt states, is a monotonically decreasing function of the mass flux through the shock, reaching an asymptotic value given by $(\Gamma-1)q\rho_0$.
On the contrary, in the relativistic regime, the Zel'dovich pressure jump shows a decreasing trend only for small values of the  mass flux. It then reaches a minimum value, followed by an increasing trend (see Fig.~\ref{fig:Pjump}). The potential impact of this new relativistic effect on the energy extraction from relativistic detonation waves has only been partially covered here, and it may deserve a separate investigation.
The present work relies on the ideal-gas equation of state both in the primitive-variable recovery of Sect.~\ref{sec:cons2prim} and in the quartic \eqref{eq:rho4} for the density that follows from the Rankine-Hugoniot conditions. More general admissible-state characterizations have also been developed for relativistic hydrodynamics with a general equation of state \citep{Wu_2017}. Moreover, robust and provably convergent primitive-variable recovery methods have recently been proposed for Synge-type equations of state \citep{YUE2026115155}. These results suggest a possible extension of the present approach beyond the ideal-gas equation of state. This analysis is beyond the scope of the present work and is left for future investigation.

Regarding the possible extension of our approach to multidimensions, 
it is easy to guess that qualitatively new effects may emerge, due to the nonlinear coupling of tangential velocities to the shock front, by means of the Lorentz factor. Preliminary results about these effects have been shown in \cite{Harpole2019}. 
Finally, we should also mention that studying reaction flows within a single fluid approximation can be regarded as a preliminary step towards a more realistic modeling where two-phase flows are considered.


\section{Acknowledgments}

We would like to thank Ilya Peshkov for inspiring discussions and two anonymous referees for their suggestions and constructive comments.
The authors of this work  are all members of the INdAM GNCS group in Italy.  

M.D. was funded by the Fondazione Caritro via the project SOPHOS and by the European Research Council (ERC) under the European Union's Horizon Europe research and innovation programme via the project SOPHOS, grant agreement no. ERC-ADG-2025-101265878-SOPHOS. 
Views and opinions expressed are however those of the authors only and do not necessarily reflect those of the European Union or the European Research Council Executive Agency. Neither the European Union nor the granting authority can be held responsible for them.

\noindent
We acknowledge the CINECA award under the ISCRA initiative, for the availability of high-performance computing (HPC) resources and technical support.

%


\bibliographystyle{plain}
\bibliography{./references.bib}


\end{document}